\documentclass[a4paper,10pt]{article}
\usepackage{amssymb}
\usepackage{latexsym}
\usepackage{amsmath}
\usepackage{amsthm}
\usepackage{amsfonts}
\usepackage{amssymb}
\usepackage{graphicx}

\newtheorem{Th}{Theorem}
\newtheorem{Prop}{Proposition}

\newtheorem{Lma}{Lemma}[section]
\newtheorem{Dfi}{Definition}

\newcommand{\be}{\begin{equation}}
\newcommand{\ee}{\end{equation}}
\newcommand{\R}{\mathbb{R}}
\newcommand{\N}{\mathbb{N}}
\newcommand{\C}{\mathbb{C}}

\newcommand\res{\mathop{\hbox{\vrule height 7pt width .5pt depth 0pt
\vrule height .5pt width 6pt depth 0pt}}\nolimits}

\newcommand{\reset}{\setcounter{equation}{0}\setcounter{Th}{0}\setcounter{Prop}{0}\setcounter{Co}{0}
\setcounter{Lm}{0}\setcounter{Rm}{0}}

\def\ti{\tilde}
\def\lf{\left}
\def\rg{\right}

\def\al{\alpha}
\def\la{\lambda}

\def\ep{\varepsilon}
\def\ds{\displaystyle}
\def\ov{\overline}

\def\om{\omega}
\def\p{\partial}

\def\res{\mathop{\hbox{\vrule height 7pt width .5pt 
depth 0pt\vrule height .5pt width 6pt depth 0pt}}\nolimits}

\begin{document}
\title{Critical weak immersed surfaces within sub-manifolds of the Teichm\"uller space.}
\author{ Tristan Rivi\`ere\footnote{Forschungsinstitut f\"ur Mathematik, ETH Zentrum,
CH-8093 Z\"urich, Switzerland.}}
%\date{ }
\maketitle

{\bf Abstract }: {\it  We prove that the critical points of various energies such as the area, the Willmore energy, the frame energy for tori...etc among possibly branched immersions constrained to evolve within a smooth sub-manifold of the Teichm\"uller space 
satisfy the corresponding constrained Euler Lagrange equation. We deduce that critical points of the Willmore energy or the frame energy for tori are smooth analytic surfaces, away possibly from isolated branched points, under the condition that either the genus is at most 2
or if the sub-manifold does not intersect  the subspace of hyper-elliptic points. Using a compactness result from \cite{Ri3} we can conclude that each closed sub-manifold of the Teichm\"uller space, including points,
under the previous assumptions, posses a possibly branched smooth Willmore minimizer satisfying the conformal-constrained Willmore equation.}

\medskip

\noindent{\bf Math. Class. 49Q10, 53A05, 53A30, 35J35}
\section{Introduction}
The purpose of the present work is to derive constrained Euler Lagrange equations for weak immersions which are critical
points of geometric energies such as the area, the Willmore energy, the frame energy for tori...etc under the
constraint that the metrics defined by the variation of this immersions stay within a given sub-manifold of the Teichm\"uller Space.\\

The study of the variations of Willmore energy has been initiated by Leon Simon in a seminal work \cite{Si} in which he was proving the existence of an embedded torus into ${\R}^m$ minimizing the $L^2$ norm of the second fundamental form. This problem was an analytical challenge at the time in particular due to the fact that this norm clearly does not provided any control the $C^1$ norm of the surface which is needed in order to speak about immersion. The need of weakening the strict geometric notion of immersion was answered in this work by considering ''measure theoretic version'' of sub-manifold known as {\it varifold} and by using local approximation of these weak objects
by  {\it bi-harmonic graphs}. In the following decade, after this analytical breakthrough, a series of works (\cite{BK}, \cite{KS1}, \cite{KS2}, \cite{KS3}...)     took over successfully this approach to solve important questions related to Willmore surfaces such as the existence of the flow, point removability property for Willmore surfaces, existence of Willmore minimizers within a given conformal class...\\

 In \cite{Ri2}, \cite{Ri2}  the author introduced a {\it parametric approach} to the study of Willmore lagrangian and the notion of {\it weak immersions}. In these works one study surfaces in ${\R}^m$ from the point of view of the immersion, the map which is generating them, while the Leon Simon's approach is mostly considering the immersed surface as a subset of the ambient space ${\R^m}$ and
is called the {\it ambiant approach} to the Willmore problem.\\

We now recall the notion of {\it weak immersion} with $L^2-$bounded second fundamental form introduced in \cite{Ri3}.\\
A closed abstract surface $\Sigma$ being given, we observe that any smooth riemannian metric is equivalent to any other one. One can
then define the Sobolev Spaces $W^{k,p}_{g_0}(\Sigma,{\R}^m)$ from $\Sigma$ into ${\R}^m$, for $k\in{\N}$ and $p\in[1,\infty]$  with respect to any smooth reference metric $g_0$ and all these spaces
are independent of  the chosen metric $g_0$ and we can simply denote them  $W^{k,p}(\Sigma,{\R}^m)$.\\

 A {\it weak immersion} 
of $\Sigma$ into ${\R}^m$ is a map
$\vec{\Phi}$ from $\Sigma$ into ${\R}^m$ such that
\begin{itemize}
\item[i)]
\[
\vec{\Phi}\in W^{1,\infty}(\Sigma,{\R}^m)\quad,
\]
\item[ii)] there exists a constant $C_{\vec{\Phi}}>1$ such that
\[
\forall x\in \Sigma\quad\forall X\in T_x\Sigma\quad\quad C_{\vec{\Phi}}^{-1}\ g_0(X)\le |d\vec{\Phi}(X)|^2\le  C_{\vec{\Phi}}\ g_0(X)
\]
i.e. in other words the metric on $\Sigma$ equal to the pull back by $\vec{\Phi}$ of the canonical metric of ${\R}^m$ is equivalent
to any reference metric on $\Sigma$,
\item[iii)]
\[
\vec{n}_{\vec{\Phi}}\in W^{1,2}(\Sigma,\wedge^{m-2}{\R}^m)
\]
where $\vec{n}_{\vec{\Phi}}$ is the Gauss map associated to $\vec{\Phi}$ and given in any local chart $(x_1,x_2)$ by
\[
\vec{n}_{\vec{\Phi}}:=\star_{{\R}^m}\frac{\p_{x_1}\vec{\Phi}\wedge\p_{x_2}\vec{\Phi}}{|\p_{x_1}\vec{\Phi}\wedge\p_{x_2}\vec{\Phi}|}
\]
where $\star_{{\R}^m}$ is the hodge operator on multi-vector associated to the canonical scalar product on ${\R}^m$.
\end{itemize}
 
The space of weak immersion is denoted ${\mathcal E}_\Sigma$ and a result whose proof goes back to the works of 
S.M\"uller and V.\v{S}ver\'ak,  \cite{MS}, and F.H\'elein, \cite{He}, asserts that any map in ${\mathcal E}_\Sigma$ defines a unique conformal class.
Precisely we have
\begin{Prop}(see \cite{Ri3} and \cite{Ri1})
\label{pr-00-1}
Let $\vec{\Phi}$ be a weak immersion in ${\mathcal E}_\Sigma$, then there exists a smooth constant scalar curvature metric $h$
of volume 1 on $\Sigma$ and a bi-lipschitz homeomorphism $\Psi$ of $\Sigma$ such that
\[
\vec{\Phi}\circ\Psi\ :\ (\Sigma,h)\ \longrightarrow\ {\R}^m\quad\quad\mbox{ is weakly conformal}\quad,
\]
more precisely there exists $\al\in L^\infty(\Sigma,{\R})\cap W^{1,2}(\Sigma,{\R})$ such that
\[
(\vec{\Phi}\circ\Psi)^\ast g_{{\R}^m}=\ e^{2\al}\ h\quad\quad\mbox{ almost everywhere}\quad.
\]
The space ${\mathcal E}_\Sigma$ defines a Banach manifold modeled on $W^{1,\infty}\cap W^{2,2}(\Sigma,{\R}^m)$. A 
family of generators of $\pi_1(\Sigma)$ being fixed, the mapping
\[
\tau\ :\ {\mathcal E}_\Sigma\ \longrightarrow\ {\mathcal T}_\Sigma
\]
which to every $\vec{\Phi}$ in ${\mathcal E}_\Sigma$ assigns the corresponding Teichm\"uller class in the Teichm\"uller Space
${\mathcal T}_\Sigma$ is smooth.\hfill $\Box$
\end{Prop}
The interest of introducing the class of weak immersions ${\mathcal E}_\Sigma$ is motivated by the following {\it almost closure
theorem}.
\begin{Th}
\cite{Ri3}
Let $\vec{\Phi}_k$ be a sequence of weak immersions in ${\mathcal E}_\Sigma$ such that
\[
\limsup_{k\rightarrow +\infty}\int_{\Sigma}|d\vec{n}_{\vec{\Phi}_k}|^2_{g_{\vec{\Phi}_k}}\ dvol_{g_{\vec{\Phi}_k}}<+\infty
\]
Assume that $\tau(\vec{\Phi}_k)$ converges in the Teichm\"uller Space ${\mathcal T}_\Sigma$ (i.e. the associated
constant scalar curvature metric $h_k$ of volume 1  converges to a limiting metric $h_\infty$) then there exists a subsequence
$\vec{\Phi}_{k'}$, a subsequence of bi-lipschitz homeomorphism $\Psi_{k'}$, a sequence $\Xi_{k'}$ of conformal transformations
of ${\R}^m\cup\{\infty\}$ and a finite number of points $\{a_1\cdots a_N\}$ such that
\[
\vec{\xi}_{k'}:=\Xi_{k'}\circ\vec{\Phi}_{k'}\circ\Psi_{k'}\ \rightharpoonup\ \vec{\xi}_\infty\quad\quad\mbox{weakly in }W^{2,2}(\Sigma\setminus\{a_1\cdots a_N\},{\R}^m)\quad
\]
where $\vec{\xi}_{k'}$ (resp. $\vec{\xi}_\infty$) is weakly conformal from $(\Sigma,h_{k'})$ (resp. $(\Sigma,h_\infty)$) into ${\R}^m$. Moreover
\[
\limsup_{k'\rightarrow+\infty}\|\log|d\xi_{k'}|_{h_{k'}}\|_{L^\infty_{loc}(\Sigma\setminus\{a_1\cdots a_N\})}<+\infty\quad,
\]
and
\[
\liminf_{k\rightarrow +\infty}\int_{\Sigma}|d\vec{n}_{\vec{\Phi}_k}|^2_{g_{\vec{\Phi}_k}}\ dvol_{g_{\vec{\Phi}_k}}\ge\int_{\Sigma}|d\vec{n}_{\vec{\xi}_\infty}|^2_{g_{\vec{\xi}_\infty}}\ dvol_{g_{\vec{\xi}_\infty}}
\]
\hfill $\Box$
\end{Th}
Let $N$ be a sub-manifold to the Teichm\"uller Space .\\

We denote by ${\mathcal E}^N_\Sigma$ the subspace of weak immersions $\vec{\Phi}$ in ${\mathcal E}_\Sigma$ such that $\tau(\vec{\Phi})\in N$.\\

The goal of the present paper is to study within ${\mathcal E}^N_\Sigma$ the variations of Lagrangians such as
\begin{itemize}
\item[i)] {\it The area}
\[
A(\vec{\Phi})=\int_\Sigma dvol_{g_{\vec{\Phi}}}
\]
\item[ii)] {\it The Willmore energy}
\[
W(\vec{\Phi})=\int_\Sigma|\vec{H}_{\vec{\Phi}}|^2\ dvol_{g_{\vec{\Phi}}}=\frac{1}{4}{\mathbb I}(\vec{\Phi})+\pi\,\chi(\Sigma)
\]
where $\vec{H}_{\vec{\Phi}}$ is the mean curvature vector of the weak immersion $\vec{\Phi}$, $\chi(\vec{\Phi})$ is the Euler characteristic of $\Sigma$ and
\[
{\mathbb I}(\vec{\Phi})=\int_{\Sigma}|d\vec{n}_{\vec{\Phi}}|^2_{g_{\vec{\Phi}}}\ dvol_{g_{\vec{\Phi}}}\quad.
\]
\item[iii)] {\it the frame energy for tori}
\[
F(\vec{\Phi},\vec{e}):=\int_{T^2}|d\vec{e}\,|_{g_{\vec{\Phi}}}^2\ dvol_{g_{\vec{\Phi}}}
\]
where $\Sigma$ is the torus $T^2$ and $\vec{e}=(\vec{e}_1,\vec{e}_2)$ is a frame of orthonormal vectors on the
tangent bundle $\vec{\Phi}_\ast( T T^2)$ : in other words $n_{\vec{\Phi}}:=\star_{{R}^m}\vec{e}_1\wedge\vec{e}_2$ where $\vec{e}_i\in S^{m-1}$ and $\vec{e}_1\cdot\vec{e}_2=0$.
\end{itemize}

As explained in \cite{MR3}, a pair $\gamma=(\gamma_1,\gamma_2)$ of two generators of the $\pi_1$ of $\Sigma$ being given and an immersion $\vec{\Phi}$
being also fixed, there exists a unique frame $\vec{e}$ critical point of $F(\vec{\Phi},\vec{e})$ such that the degree of the frame
along the generators $\gamma$ is zero and it is minimizing $F(\vec{\Phi},\vec{e})$ in this class. We denote by $F_\gamma(\vec{\Phi})$
the $F$ energy for this particular frame.\\

These energy are just examples in order to illustrate our main result. In order to state this main result we need to introduce the notions of {\it Weingarten form}, {\it holomorphic quadratic differential} and {\it Weil Peterson hermitian product} for a weak immersion $\vec{\Phi}$ of a 2-manifold $\Sigma$.\\

 Considering a smooth immersion $\vec{\Phi}$ of an arbitrary 2-dimensional manifold $\Sigma$
into ${\R}^m$ one can then introduce the corresponding bundle of $1-0$ forms over $\Sigma$ denoted $\wedge^{1-0}T^\ast\Sigma$.\\

One defines the {\it Weingarten form} using local conformal charts given by  proposition~\ref{pr-00-1}
as being the following global section of ${\R}^m\otimes\wedge^{1-0}T^\ast\Sigma\otimes\wedge^{1-0}T^\ast\Sigma$  :
\be
\label{0I.9}
\begin{array}{l}
\vec{h}_0:= 2\,e^{-2\la}\,\pi_{\vec{n}}(\p^2_{z^2}\vec{\Phi})\ dz\otimes dz\\[5mm]
\ds\quad=\frac{e^{-2\la}}{2}\,\pi_{\vec{n}}\lf(\p^2_{x^2_1}\vec{\Phi}-\p^2_{x^2_2}\vec{\Phi}-2\, i\ \p^2_{x_1x_2}\vec{\Phi}\rg)\ \ {dz\otimes dz}
\end{array}
\ee
where $\pi_{\vec{n}}$ is the orthogonal projection onto the plane orthogonal to $\vec{\Phi}_\ast T\Sigma$.\\

Holomorphic quadratic forms associated to $\vec{\Phi}$ are holomorphic sections of  $\wedge^{1-0}T^\ast\Sigma\otimes\wedge^{1-0}T^\ast$.\\

The class $\tau(\vec{\Phi})$ in the Teichm\"uller space ${\mathcal T}_\Sigma$ of $\Sigma$ associated to the immersion $\vec{\Phi}$
is said to be {\it hyperelliptic} if the tensor products of holomorphic 1-forms do not generate the vector space of holomorphic quadratic
form. We recall that for $g\ge 3$ the subset $H_g$ of hyper-elliptic classes is a complex analytic sub-manifold of ${\mathcal T}_\Sigma$ 
of complex co-dimension $g-2$  (see \cite{Na} 4.1.5).\\
Let $\vec{\Phi}$ be a weak immersion of $\Sigma$. We express this immersion in a conformal chart from the 2-disc  $D^2$
and we keep denoting $\vec{\Phi}$ this mapping.  We introduce on the space $\wedge^{1-0}D^2\otimes\wedge^{1-0} D^2$ of $1-0\otimes1-0$ form on $D^2$ the following hermitian product depending on the conformal immersion $\vec{\Phi}$
\[
(\psi_1\ dz\otimes dz,\psi_2\ dz\otimes dz)_{WP}:= e^{-4\la}\ \ov{\psi_1(z)}\ \psi_2(z)
\]
where $e^\la:=|\p_{x_1}\vec{\Phi}|=|\p_{x_2}\vec{\Phi}|$. We observe that for a conformal change of coordinate $w(z)$ (i.e. $w$ is holomorphic in $z$) and for $\psi_i'$ satisfying
\[
\psi_i'\circ w\ dw\otimes dw=\psi_i\ dz\otimes dz
\]
one has, using the conformal immersion $\vec{\Phi}\circ w$ in the l.h.s.
\[
(\psi_1'\ dw\otimes dw,\psi_2'\ dw\otimes dw)_{WP}=(\psi_1\ dz\otimes dz,\psi_2\ dz\otimes dz)_{WP} 
\]
Hence this hermitian product is independent of the conformal chart and only depends on the conformal structure defined by $\vec{\Phi}$.
The scalar product given by the imaginary part of $<\cdot,\cdot>_{WP}$ is denoted
\[
<\cdot,\cdot>_{WP}:=\Im(\cdot,\cdot)_{WP}
\]
Integrated on $\Sigma$ it defines  the so called {\it Weil Peterson product}.

\begin{Th}
\label{th-0-1}
Let $\Sigma$ be a closed two dimensional manifold. Let $N$ be a sub-manifold of the Teichm\"uller space ${\mathcal T}_\Sigma$. Let $\tau_0\in {\mathcal T}_\Sigma$
and assume that either $g\le 2$  or $\tau_0$ is not hyper-elliptic. Let $\vec{\Phi}$ be a weak immersion in ${\mathcal E}_\Sigma$
from $\Sigma$ into ${\R}^m$
such that $\tau(\vec{\Phi})=\tau_0$. Assume $\vec{\Phi}$ is a critical 
point  of
the area, resp. the Willmore energy, resp. the frame energy for all $C^1$ perturbations included in $N$ then there exists an holomorphic quadratic form $q$ of $(\Sigma,\tau_0)$ (i.e. an holomorphic section of
$\wedge^{1,0}T^\ast\Sigma\otimes \wedge^{1,0}T^\ast\Sigma$) such that
\begin{itemize} 
\item[i)]
\be
\label{00-1}
\vec{H}=<q,\vec{h}_0>_{WP}
\ee
\item[ii)]
\be
\label{00-2}
d^{\ast_g}\lf[ d\vec{H}-3\, \pi_{\vec{n}}(d\vec{H})+\star_{{\R}^m}(\ast_{g}d\vec{n}_{\vec{\Phi}}\wedge\vec{H})\rg]=<q,\vec{h}_0>_{WP}
\ee
where $\ast_g$ is the Hodge operator on $\wedge^pT^\ast\Sigma$ issued by $g$.
\item[iii)]
\be
\label{00-3}
\lf\{
\begin{array}{l}
\ds d^{\ast_g}\lf[ d\vec{H}-3\, \pi_{\vec{n}}(d\vec{H})+\star_{{\R}^m}(\ast_{g}d\vec{n}_{\vec{\Phi}}\wedge\vec{H})
-\vec{\mathbb I}\res_g(\vec{e}_2\cdot d\vec{e}_1)\rg.\\[5mm]
\ds\quad\lf.+[\vec{e}_2\cdot d\vec{e}_1\otimes\vec{e}_2\cdot d\vec{e}_1-2^{-1}|\vec{e}_2\cdot d\vec{e}_1|^2]\res_gd\vec{\Phi} \rg]=<q,\vec{h}_0>_{WP}\\[5mm]
\ds\ d^\ast(\vec{e}_2\cdot d \vec{e}_1)=0
\end{array}
\rg.
\ee
where $\res_g$ is the standard contraction operator between  a $p-$vectors and a $q-$vectors  $(p\ge q)$ given by
\[
\forall\,\vec{a}\in\wedge^p{\R}^m\ ,\quad\forall\,\vec{b}\in\wedge^q{\R}^m\, ,\quad\forall\,\vec{c}\in\wedge^{p-q}{\R}^m\quad\quad<\vec{a}\res\vec{b},\vec{c}>_g=<\vec{a},\vec{b}\wedge\vec{c}>_g\quad.
\] 

\end{itemize}
\hfill $\Box$
\end{Th}
This result has been established when $\vec{\Phi}$ is a non-degenerate critical point of $\tau$ and when $N$ is a point
 in \cite{BPP} for smooth immersions and in \cite{Ri3} for weak immersions. Being a degenerate critical point
of the constraint $\tau(\vec{\Phi})=\tau_0$ is equivalent to the fact that there exists a non-zero holomorphic quadratic
form $q$ such that
\[
<q,\vec{h}_0>_{WP}\equiv 0\quad\quad\mbox{ in }{\mathcal D}'(\Sigma)\quad.
\]
Such a weak immersion is called {\it global isothermic} (see \cite{Ri5}). In \cite{KS4}, using quite involved calculations, the authors have been able to treat the degenerate case when $\vec{\Phi}$ is a smooth isothermic immersion when the constraint
is reduced to a point in the Teichm\"uller space.
The main achievement of the present work is to present a general argument for dealing with perturbations of $C^1$ energies
at arbitrary  global isothermic surfaces in arbitrary co-dimension and within an arbitrary sub-manifold of the Teichm\"uller space avoiding the hyper-elliptic points in case when the genus
of $\Sigma$ is larger than 2.\\

Using the regularity theory for weak immersions satisfying (\ref{00-2}) and (\ref{00-3}) respectively in \cite{Ri2}, \cite{Ri3} and \cite{MR3}
we deduce the following
\begin{Th}
\label{th-0-2}
Let $\Sigma$ be a closed two dimensional manifold. Let $N$ be a sub-manifold of the Teichm\"uller space ${\mathcal T}_\Sigma$. Let $\tau_0\in {\mathcal T}_\Sigma$
and assume that either $g\le 2$  or $\tau_0$ is not hyper-elliptic. Let $\vec{\Phi}$ be a weak immersion in ${\mathcal E}_\Sigma$
from $\Sigma$ into ${\R}^m$
such that $\tau(\vec{\Phi})=\tau_0$. Assume $\vec{\Phi}$ is a critical 
point  of the Willmore energy, resp. the frame energy for all $C^1$ perturbations included in $N$ then $\vec{\Phi}$ is analytic.
\hfill $\Box$
\end{Th}

For variational purposes it is convenient to extend a bit the class of weak immersion allowing isolated branched points singularities. For a given riemann surface $(\Sigma,h)$ equipped
with a compatible constant gauss curvature metric $h$ we define respectively
\[
{\mathcal F}_{(\Sigma,h)}^{conf}:=\lf\{
\begin{array}{l}
\ds\vec{\Phi}\in W^{1,\infty}(\Sigma,{\R}^m)\quad\mbox{s.t. }\exists \,b_1\cdots b_n\in \Sigma\\[5mm]
\ds\quad\exists\, \al\in L^\infty_{loc}(\Sigma\setminus\{b_1\cdots b_n\})\quad,\quad \vec{\Phi}^\ast g_{{\R}^m}=e^{2\al}\, h\\[5mm]
\ds\quad e^{2\al}\in L^1(\Sigma,h)\quad \mbox{ and }\quad d\vec{n}_\Phi\in L^2(\Sigma,h)
\end{array}
\rg\}
\]
and the space of weak branched immersion ${\mathcal F}_\Sigma$
\[
{\mathcal F}_\Sigma=\lf\{
\begin{array}{l}
\vec{\Phi} \in W^{1,\infty}(\Sigma,{\R}^m)\quad\mbox{s.t. }\quad\exists\, \Psi \mbox{ bi-lipschitz homeo. of }\Sigma\\[5mm]
\mbox{ and } \exists\,h\mbox{ metric of constant curvature s.t. }\quad \vec{\Phi}\circ\Psi\in {\mathcal F}_{(\Sigma,h)}^{conf}
\end{array}
\rg\}
\]
The following theorem holds
\begin{Th}
\label{th-0-3} (\cite{Ri2},\cite{MR2} see also \cite{Ri1}) ${\mathcal F}_\Sigma$ is weakly sequentially complete under uniform
$L^2$ control of the second fundamental form and control of the Teichm\"uller class. In other words, 
let $\vec{\Phi}_k$ be a sequence of branched weak immersions in ${\mathcal F}_\Sigma$ such that
\[
\limsup_{k\rightarrow +\infty}\int_{\Sigma}|d\vec{n}_{\vec{\Phi}_k}|^2_{g_{\vec{\Phi}_k}}\ dvol_{g_{\vec{\Phi}_k}}<+\infty
\]
Assume that $\tau(\vec{\Phi}_k)$ converges in the Teichm\"uller Space ${\mathcal T}_\Sigma$ (i.e. the associated
constant scalar curvature metric $h_k$ of volume 1  converges to a limiting metric $h_\infty$) then there exists a subsequence
$\vec{\Phi}_{k'}$, a subsequence of bi-lipschitz homeomorphism $\Psi_{k'}$, a sequence $\Xi_{k'}$ of conformal transformations
of ${\R}^m\cup\{\infty\}$ and a finite number of points $\{a_1\cdots a_N\}$ such that
\[
\vec{\xi}_{k'}:=\Xi_{k'}\circ\vec{\Phi}_{k'}\circ\Psi_{k'}\ \rightharpoonup\ \vec{\xi}_\infty\quad\quad\mbox{weakly in }W^{2,2}(\Sigma\setminus\{a_1\cdots a_N\},{\R}^m)\quad
\]
where $\vec{\xi}_{k'}\in{\mathcal F}^{conf}_{(\Sigma,h_\infty)}$. Moreover
\[
\limsup_{k'\rightarrow+\infty}\|\log|d\xi_{k'}|_{h_{k'}}\|_{L^\infty_{loc}(\Sigma\setminus\{a_1\cdots a_N\})}<+\infty\quad,
\]
and
\[
\liminf_{k\rightarrow +\infty}\int_{\Sigma}|d\vec{n}_{\vec{\Phi}_k}|^2_{g_{\vec{\Phi}_k}}\ dvol_{g_{\vec{\Phi}_k}}\ge\int_{\Sigma}|d\vec{n}_{\vec{\xi}_\infty}|^2_{g_{\vec{\xi}_\infty}}\ dvol_{g_{\vec{\xi}_\infty}}\quad.
\]

\hfill $\Box$
\end{Th}
The map $\tau$ into the Teichm\"uller space ${\mathcal T}_\Sigma$ is naturally extended to ${\mathcal F}_\Sigma$.
Similarly as in the case of weak immersions in ${\mathcal E}_\Sigma$,  a sub-manifold $N$ of the Teichm\"uller space ${\mathcal T}_\Sigma$ being given, we denote
by ${\mathcal F}^N_{\Sigma}$ the subspace of ${\mathcal F}_\Sigma$ made of branched weak immersions $\vec{\Phi}$ such that $\tau(\vec{\Phi})\in N$.
Combining now the previous completeness result, the regularity result theorem~\ref{th-0-2}, together with the Frechet differentiability of ${\mathbb I}$
and $F$ proved in \cite{Ri3} and \cite{MR2}, we obtain the following theorem which one of the consequences of the main result, theorem~\ref{th-0-1} of the present work.
\begin{Th}
\label{th-0-4}
Let $\Sigma$ be a closed two dimensional manifold. Let $N$ be a closed (compact without boundary) sub-manifold of the Teichm\"uller space ${\mathcal T}_\Sigma$. 
Then there is an absolute minimizer of ${\mathbb I}$ within $N$ in the space
of branched weak immersions ${\mathcal F}^N_\Sigma$,  it satisfies the Willmore conformally constrained equation (\ref{00-2})
and is analytic away from possibly finitely many branched points. Similarly, taking two generators $\gamma=(\gamma_1,\gamma_2)$ of  $\pi_1(\Sigma)$, The frame energy for frames $\vec{e}$ which have zero degrees along these curves admits an analytic minimizing immersion  in ${\mathcal E}^N_\Sigma$ and it satisfies the frame constrained equation (\ref{00-3})
\hfill$\Box$
\end{Th}
Existing result of minimizers of ${\mathbb I}$  within a conformal class has been obtained in \cite{Ri3}. Here it has been proved that either the minimizer is analytic  away from possibly finitely many branched points and satisfies the Willmore conformally constrained equation (\ref{00-2}) or is {\it isothermic}. Independently, the existence of a minimizer
in the space of weak immersions within a conformal class, under some energy bound assumption, has been obtained in \cite{KL}. In \cite{KS4} the authors have been able to prove, in codimension 1 and 2, under some energy bound assumption, that the minimizer is always satisfying the {\it conformal constrained} equation (\ref{00-2}) away from branched points
(recall the result from \cite{Ric} saying that, in 3-dimension, {\it isothermic} smooth immersions satisfying the {\it conformal constrained} equation (\ref{00-2}) are {\it constant mean curvature}
surfaces in a space form).

\section{Notations - Preliminary results.}
\reset
\subsection{The Period Matrix $\Pi(h)$.}
Let $\Sigma$ be a closed 2-dimensional orientable manifold. Let $g$ be the genus of $\Sigma$. We assume $g\ge 1$.

\medskip

Let $a_1\cdots a_g\, , b_1\cdots b_g$ be a canonical basis for the homology $H^1(\Sigma)$ : it satisfies for any $j,k=1\cdots g$
\be
\label{0.0}
\lf\{
\begin{array}{l}
a_j\cdot b_k=\delta_{jk}\\[5mm]
a_j\cdot a_k=b_j\cdot b_k=0\quad.
\end{array}
\rg.
\ee
Let $h_0$ be a metric on $\Sigma$. Denote by $\al^k_0$ a family of harmonic 1-forms on $(\Sigma,h_0)$ (see proposition III.2.8 page 63 of \cite{FK}) such that
\be
\label{0.1}
\int_{a_j}\al^k_0=\delta_{jk}\quad\mbox{ and }\int_{a_j}\ast_{0}\,\al^k_0=0\quad,
\ee
where $\ast_0$ is the Hodge operator associated to $h_0$ and such that
\be
\label{0.2}
\lf(\int_{b_j}\al^k_0\rg)_{j,k=1\cdots g}\quad\quad\mbox{ is symmetric}
\ee
and
\be
\label{0.3}
\lf(\int_{b_j}\ast_0\al^k_0\rg)_{j,k=1\cdots g}\quad\quad\mbox{ is symmetric positive definite}\quad.
\ee
We consider metrics $h$ in the $L^\infty$ neighborhood of $h_0$ and we denote 
\be
\label{II.1}
\al^k(h):=\al^k_0+d\varphi^k(h)
\ee
where $\varphi^k(g)$ is the solution to
\be
\label{II.1a}
\Delta_{h}\varphi^k(h)=-d^{\ast_h}\al^k_0\quad.
\ee
and $\Delta_h$ is the Laplace Beltrami operator on $(\Sigma,h)$ and $d^{\ast_h}$ is the adjoint operator to $d$ for the metric $h$. Denote by $\ast_h$ the 
Hodge operator associated to $h$ - we have then $\ast_0=\ast_{h_0}$ - and by 
\[
\om^k(h):=\al^k(h)+i\ast_h\al^k(h)\quad\quad k=1\cdots g
\]
the basis of holomorphic 1-forms (abelian differential of first kind) associated to $\al^k(h)$. We denote by $\Pi(h)$ the period map
\[
\Pi(h):=(\Pi^0(h),\Pi^1(h))=\lf(\lf(\int_{a_j}\om^k(h)\rg)_{j,k=1\cdots g},\lf(\int_{b_j}\om^k(h)\rg)_{j,k=1\cdots g}\rg)
\]
and
\[
\Pi_{jk}^0(h):=\int_{a_j}\om^k(h)\quad\quad\mbox{ and }\quad\quad\Pi_{jk}^1(h):=\int_{b_j}\om^k(h)
\]
So that we have
\[
\Pi_{jk}^0(h_0)=\delta_{jk}\quad\quad\mbox{ and }\quad\quad\Pi_{jk}^1(h_0)= c_{jk}+i\,d_{jk}
\]
where $D:=(d_{jk})_{j,k=1\cdots g}$ is invertible. Observe that the choice we are making of the variation of basis of holomorphic 1-form around $h_0$ imposes
\be
\label{II.2}
\Re(\Pi_{jk}^0(h))\equiv \delta_{jk}\quad\quad\mbox{ and }\quad\quad\Re(\Pi_{jk}^1(h))\equiv c_{jk}\quad.
\ee
We also denote by $\tau(h)$ the {\it Teichm\"uller class} induced by $h$ for the choice of basis we have made in $H_1(\Sigma)$.

\medskip

A theorem by Torelli asserts the following (see for instance \cite{Na})
\begin{Th}
\label{th-II.1}
Let $\Sigma$ be a two dimensional closed orientable manifold, let $a_1\cdots a_g,\, b_1\cdots b_g$ be a canonical basis for $H_1(\Sigma)$. Let $\tau$ and $\tau'$ be two Teichm\"uller classes 
on $\Sigma$ and let $\om=(\om_1\cdots\om_g)$ and $\om':=(\om'_1\cdots\om_g')$ be the two basis of holomorphic 1-forms  - see proposition III.2.8 page 63 of \cite{FK}) -
such that the corresponding period matrix satisfying respectively
\[
\Pi=(I_g,\Pi^1)\quad\quad\mbox{ and }\quad\quad\Pi'=(I_g,(\Pi^1)')
\]
where $\Pi^1$ and $(\Pi^1)'$ are symmetric matrices with positive-definite imaginary parts, then 
\[
\Pi^1=(\Pi^1)'\quad\quad\Leftrightarrow\quad\quad \tau=\tau'\quad.
\]
\hfill $\Box$
\end{Th}
For $h$ in an $L^\infty$ neighborhood of $h_0$ the matrix $\Pi^0$ is invertible and if we choose
\[
\ti{\om}_k(h):=\sum_{j=1}^g e_{jk}(h)\,\om_j(h)\quad,
\]
where $E(g)=(e_{jk}(h))_{j,k=1\cdots g}=((\Pi)^0)^{-1}$, the corresponding period function for the canonical basis $a_1\cdots a_g,\, b_1\cdots b_g$ is given by
\[
\ti{\Pi}(h)= (Id_g,\ti{\Pi}^1(h))=(Id_g, E(h)\,\Pi^1(h))
\]
For two metrics $h$ and $h'$ in the neighborhood of $h_0$ we have
\[
\Pi(h)=\Pi(h')\quad\Rightarrow\quad E(h)\,\Pi^1(h)=E(h')\,\Pi^1(h')\quad\Rightarrow\quad \ti{\Pi}^1(h)=\ti{\Pi}^1(h')\quad\Rightarrow\quad \tau(h)=\tau(h')\quad.
\]
 We assume now that $\tau(h)=\tau(h')$.
 
 \medskip
 
  Let $\eta_{a_j}(h)$ and $\eta_{b_j}(h)$ be the harmonic 1-form for $h$ representing the Poincar\'e dual of respectively $a_j$ and $b_j$. We have in particular
 (see for instance \cite{FK} chapter III)
 \[
 \lf\{
 \begin{array}{l}
\ds \int_\Sigma \eta_{a_j}(h)\wedge\eta_{b_k}(h)=\int_{a_j}\eta_{b_k}(h)=-\int_{b_k}\eta_{a_j}(h)=a_j\cdot b_k=\delta_{jk}\\[5mm]
\ds\int_\Sigma \eta_{a_j}(h)\wedge\eta_{a_k}(h)=\int_{a_j}\eta_{a_k}(h)=-\int_{a_k}\eta_{a_j}(h)=a_j\cdot a_k=0\\[5mm]
\ds\int_\Sigma \eta_{b_j}(h)\wedge\eta_{b_k}(h)=\int_{b_j}\eta_{b_k}(h)=-\int_{b_k}\eta_{b_j}(h)=b_j\cdot b_k=0
\end{array}
\rg.
 \]
Taking $\beta^k(h):=\eta_{b_k}(h)+\sum_{l=1}^gc_{lk}\,\eta_{a_l}(h)$  we have
\be
\label{II.3}
\forall\, k,j=1\cdots g\quad\quad,\quad\int_{a_j}\al^k(h)-\beta^k(h)=0\quad\mbox{ and }\int_{b_j}\al^k(h)-\beta^k(h)=0
\ee
Thus $\al^k(h)-\beta^k(h)$ are zero cohomologic and harmonic, this implies that $\al^k(h)=\beta^k(h)$ and in particular
\be
\label{II.3a}
\om^k(h)=\eta_{b_k}(h)+i\ast_h\,\eta_{b_k}(h)+\sum_{l=1}^gc_{lk}\,[\eta_{a_l}(h)+i\ast_h\,\eta_{a_l}(h)]
\ee
Since $\tau(h)=\tau(h')$, there exists a conformal diffeomorphism homotopic to the identity $\psi$ from $(\Sigma,h')$ into $(\Sigma,h)$.
Since $\Psi$ is homotopic to the identity, $\Psi^{-1}$ induced the identity on in $H_1(\Sigma)$ i.e. $\Psi^{-1}_\ast a_j=a_j$ and $\Psi^{-1}_\ast b_j=b_j$. Thus we have
\[
\lf\{
\begin{array}{l}
\ds\int_{a_j}\Psi^\ast\eta_{b_k}(h)= \int_{\Psi^{-1}_\ast a_j}\Psi^\ast\eta_{b_k}(h)=\int_{a_j}\eta_{b_k}(h)=\delta_{jk}\quad,\\[5mm]
\ds\int_{b_j}\Psi^\ast\eta_{a_k}(h)= \int_{\Psi^{-1}_\ast b_j}\Psi^\ast\eta_{a_k}(h)=\int_{b_j}\eta_{a_k}(h)=-\delta_{jk}\quad,\\[5mm]
\ds\int_{a_j}\Psi^\ast\eta_{a_k}(h)= \int_{\Psi^{-1}_\ast a_j}\Psi^\ast\eta_{a_k}(h)=\int_{a_j}\eta_{a_k}(h)=0\quad,\\[5mm]
\ds\int_{b_j}\Psi^\ast\eta_{b_k}(h)= \int_{\Psi^{-1}_\ast b_j}\Psi^\ast\eta_{b_k}(h)=\int_{b_j}\eta_{b_k}(h)=0\quad,\\[5mm]
\end{array}
\rg.
\]
Thus the cohomology class of the closed forms $\Psi^\ast\eta_{a_k}(h)$ (resp. $\psi^\ast\eta_{b_k}(h)$) are the Poincar\'e duals of $a_k$ (resp. $b_k$).
Since  $\eta_{a_k}(h)$ is harmonic in $(\Sigma,h)$, $\ast_h\,\eta_{a_k}(h)$ is closed and the 1 forms $\Psi^\ast\ast_h\,\eta_{a_k}(h)$ are closed as well (idem for $\Psi^\ast\ast_h\,\eta_{b_k}(h)$). Moreover since $\Psi$ is conformal we have
\[
\Psi^\ast \ast_h\eta_{a_k}(h)=\ast_{h'}\Psi^\ast \eta_{a_k}(h)\quad\quad\mbox{ and }\quad\quad\Psi^\ast \ast_h\eta_{b_k}(h)=\ast_{h'}\Psi^\ast \eta_{b_k}(h)
\] 
Thus $\ast_{h'}\Psi^\ast \eta_{a_k}(h)$ are closed (idem for $\ast_{h'}\Psi^\ast \eta_{b_k}(h)$) which implies that $\Psi^\ast \eta_{a_k}(h)$ (resp. $\Psi^\ast \eta_{b_k}(h)$) is the harmonic representative for $h'$ of the Poincar\'e dual to $a_k$ (resp. $b_k$). In other words we have proved
\be
\label{II.4}
\Psi^\ast \eta_{a_k}(h)=\eta_{a_k}(h')\quad\quad\mbox{ and }\quad\quad\Psi^\ast \eta_{b_k}(h)=\eta_{b_k}(h')\quad.
\ee
Because of (\ref{II.3a}) we have
\be
\label{II.3b}
\Psi^\ast\om^k(h)=\om^k(h')
\ee 
This implies finally that $\Pi(h)=\Pi(h')$. We have then established the following proposition.

\begin{Prop}
\label{pr-II.1}
In an $L^\infty$ neighborhood of the metric $h_0$, under the notation above, two different metrics $h$ and $h'$ define
the same Teichm\"uller class $\tau=\tau(h)=\tau(h')$ if and only if
\be
\label{II.5}
\Pi(h)=\Pi(h')\quad,
\ee
Thus $\Pi$ can be seen as a function of $\tau$ for the Teichm\"uller classes in a neighborhood of $\tau(h^0)$.
\hfill $\Box$
\end{Prop}

A classical result by Ahlfors (see \cite{Ahl}) asserts that there exists a complex structure on the Teichm\"uller space $T(\Sigma)$ - equiped with he Teichm\"uller topology -  of $\Sigma$ such
that the period map $\ti{\Pi}^1$ - viewed as a map from $T(\Sigma)$ into $Sym_{\Im,+}({\C},g)$, the space of complex  symmetric matrices with definite positive imaginary part - is holomorphic moreover if $h_0$ is not defining an {\it hyperelliptic riemann surface}\footnote{ We recall that a riemann surface is hyperelliptic if the squares of holomorphic 1-forms do not generate the 
space of holomorphic quadratic forms of the surface. For $g>1$ the subspace of hyperelliptic surfaces is an holomorphic submanifold of $T(\Sigma)$ with countably many components
- see again \cite{Ahl}.} or if $g=2,1$
one has that the complex rank to $d\ti{\Pi}^1)_{\tau(h_0)}$ is maximal and equals the complex dimension of $T(\Sigma)$
\be
\label{II.6}
rank_{{\C}}( d\ti{\Pi}^1)_{\tau(h_0)}\lf\{
\begin{array}{l}
=3 g-3\quad\quad g\ge 2\\[5mm]
=1\quad\quad g=1\quad.
\end{array}
\rg.
\ee
Whereas, if $\tau(h_0)$ corresponds to an  {\it hyperelliptic riemann surface} one has
\be
\label{II.7}
rank_{{\C}}( d\ti{\Pi}^1)_{\tau(h_0)}=2g-1\quad.
\ee
Since $\Pi=\Pi^0\,\ti{\Pi}=\Pi^0\,(Id_g,\ti{\Pi}^1)=(\Pi^0,\Pi^1)$ we have that for any variation $\tau$ in the Teichm\"uller space
\[
d_\tau\Pi(\tau(h^0))=\lf(d_\tau\Pi^0(\tau(h^0)),d_\tau\Pi^0(\tau(h^0))\,\ti{\Pi}^1(\tau(h^0))+\Pi^0(\tau(h^0))\, d_\tau\ti{\Pi}^1(\tau(h^0))\rg)
\]
Assume $\tau\in Ker d\Pi(\tau(h^0))$ then we must have 
$$
d_\tau\Pi^0(\tau(h^0))=0\quad\quad\mbox{ and }\quad\quad d_\tau\Pi^0(\tau(h^0))\,\ti{\Pi}^1(\tau(h^0))+ d_\tau\ti{\Pi}^1(\tau(h^0))=0
$$
where we have used the fact that $\Pi^0(\tau(h^0))=Id_g$. This implies that $d_\tau\ti{\Pi}^1(\tau(h^0))=0$. Thus we deduce that in the non hyperelliptic case
\be
\label{II.6a}
rank_{{\C}}( d{\Pi})_{\tau(h_0)}\lf\{
\begin{array}{l}
=3 g-3\quad\quad g\ge 2\\[5mm]
=1\quad\quad g=1\quad.
\end{array}
\rg.
\ee
Whereas, if $\tau(h_0)$ corresponds to an  {\it hyperelliptic riemann surface} one has
\be
\label{II.7a}
rank_{{\C}}( d{\Pi})_{\tau(h_0)}=2g-1\quad.
\ee
\subsection{ Computation of $d\Pi(h)$ and $d^2\Pi(h)$ at the origin $h^0$}

we shall consider only variations of the metric supported in a single chart that we choose to be complex
for the complex structure induced by $h^0$. In this chart we write
\[
h^0=e^{2\la}\ (dx_1^2+dx_2^2)\quad.
\]
and we will look at $h=h^0+ \nu$ where $\nu=\delta h$ is an arbitrary map compactly supported in the chart taking values into symmetric $2\times 2$ matrices.
 In this chart  we write
 \[
 \al^k(h)= X^k_1(h)\, dx_1+X^k_2(h)\, dx_2=(X^k_1(h^0)+\p_{x_1}\varphi^k(h))\ dx_1+  (X^k_2(h^0)+\p_{x_2}\varphi^k(h))\ dx_2               
 \] 
and
\[
\eta_{a_l}=A^l_1\ dx_1+A^l_2\ dx_2\quad\quad\mbox{ and }\eta_{b_l}=B^l_1\ dx_1+B^l_2\ dx_2\
\]
We also denote
\[
\ast_{h}\ dx_j=\sum_{i=1}^2 J^{i,j}(h)\ dx_i\quad\quad\mbox{ and }\quad\quad J(h)=(J^{i,j}(h))_{i,j=1,2}
\]
Since the coordinates are complex for $h^0$, we have
\[
J(h^0)=J_0:=\lf(
\begin{array}{cc}
0 & -1\\[5mm]
1 & 0
\end{array}\rg)
\]
With these notations we have
\[
\ast_h\ \al^k(h)=\sum_{i,j=1}^2J^{i,j}(h)\, X^k_j(h)\ dx_i
\]
and 
\[
\eta_{a_l}(h)\wedge\ast_h\ \al^k(h)= ( A^l_1\ dx_1+A^l_2\ dx_2)\wedge   \lf(\sum_{i,j=1}^2J^{i,j}(h)\, X^k_j(h)\ dx_i\rg)=(A^l)^t\, J_0^{t}\, J(h)\, X^k(h)\ dx^2
\]
where $dx^2$ denotes the canonical 2-form $dx^2=dx_1\wedge dx_2$. We have that
\be
\label{II.8}
\p_\nu \Pi^0_{lk}(h)=i\,\int_{\Sigma}\eta_{a_l}\wedge \ast_h d(\p_\nu\varphi^k(h))+i\,\int_{D^2}(A^l)^t\, J_0^{t}\, \p_\nu J(h)\, X^k(h)\ dx^2\quad,
\ee
and
\be
\label{II.9}
\p_\nu \Pi^1_{lk}(h)=i\,\int_{\Sigma}\eta_{b_l}\wedge \ast_h d(\p_\nu\varphi^k(h))+i\,\int_{D^2}(B^l)^t\, J_0^{t}\, \p_\nu J(h)\, X^k(h)\ dx^2\quad,
\ee
In particular, at $h=h_0$, since $\ast_{0}\eta_{a_i}$ and $\ast_{0}\eta_{b_i}$ are closed we have
\be
\label{II.8a}
\p_\nu \Pi^0_{lk}(h^0)=i\,\int_{D^2}(A^l)^t\, J_0^{t}\, \p_\nu J(h^0)\, X^k_0\ dx^2\quad,
\ee
and
\be
\label{II.9a}
\p_\nu \Pi^1_{lk}(h^0)=i\,\int_{D^2}(B^l)^t\, J_0^{t}\, \p_\nu J(h^0)\, X^k_0\ dx^2\quad.
\ee
For the same reason the second derivatives of $\Pi$ at $h^0$ are given by
\be
\label{II.8b}
\p^2_{\nu^2} \Pi^0_{lk}(h^0)=i\,\int_{D^2}(A^l)^t\, J_0^{t}\, \p^2_{\nu^2} J(h^0)\, X^k_0\ dx^2+2i\,\int_{D^2}(A^l)^t\, J_0^{t}\, \p_{\nu} J(h^0)\, \nabla(\p_\nu \varphi^k(h^0))\ dx^2\quad,
\ee
and
\be
\label{II.9b}
\p^2_{\nu^2} \Pi^1_{lk}(h^0)=i\,\int_{D^2}(B^l)^t\, J_0^{t}\, \p^2_{\nu^2} J(h^0)\, X^k_0\ dx^2+2i\,\int_{D^2}(B^l)^t\, J_0^{t}\, \p_{\nu} J(h^0)\, \nabla(\p_\nu \varphi^k(h^0))\ dx^2\quad.
\ee
We are now going to express $J^t_0\p_\nu J(h^0)$ and $J^t_0\p^2_{\nu^2}J(h^0)$ in terms of $h^0$ and $\nu$. Let $$G(h)=(h_{ij})_{i,j=1,2}=(e^{2\la}\, \delta_{ij}+\nu_{ij})_{i,j=1,2}$$ be the expression of $h$ in the coordinates.
The classical definition for $\ast_h$ says that for any pair of  1-forms $\al$ and $\beta$ 
\[
\al\wedge \ast_h\beta= h(\al,\beta)\ dvol_h\quad.
\]
Writing $\al=X_1\, dx_1+X_2\, dx_2$ and $\beta=Y_1\, dx_1+Y_2\, dx_2$, the previous identity becomes
\be
\label{II.10}
X^t\,J_0^t\,J(h)\, G(h)\,Y=X^t\, Y\ \sqrt{det(G(h))}\quad.
\ee
Thus we have
\be
\label{II.11}
J_0^t\,\p_{\nu}J(h)\, G(h)+J_0^t\,J(h)\,\p_\nu G(h)= I_2\ \p_\nu\sqrt{det(G(h))}\quad,
\ee
and
\be
\label{II.12}
e^{2\la}\, J_0^t\,\p^2_{\nu^2}J(h^0)+2\,J_0^t\,\p_{\nu}J(h^0)\, \nu= I_2\ \p^2_{\nu^2}\sqrt{det(G)}(h^0)\quad,
\ee
where we have used the fact that $\p^2_{\nu^2} G(h)=0$. From (\ref{II.11}) we obtain
\be
\label{II.13}
e^{2\la}\,J_0^t\,\p_{\nu}J(h^0)+\nu= I_2\ \p_\nu\sqrt{det(G)}(h^0)\quad.
\ee
We have
\be
\label{II.14}
\p_\nu\sqrt{det(G(h))}=\frac{\p_\nu det(G(h))}{2\,\sqrt{det(G(h))}}\quad,
\ee
and
\be
\label{II.15}
\p^2_{\nu^2}\sqrt{det(G(h))}=\frac{\p^2_{\nu^2} det(G(h))}{2\,\sqrt{det(G(h))}}-\frac{(\p_\nu det(G(h)))^2}{4\,(det(G(h)))^{3/2}}\quad.
\ee
For $2\times 2$ matrices one has
\be
\label{II.16}
det(G(h_0)+\nu)=det( G(h_0))+det(\nu)+trG(h_0)\, tr\nu-tr(G(h_0)\,\nu)
\ee
Combining (\ref{II.14}), (\ref{II.15}) and (\ref{II.16}) gives in one hand
\be
\label{II.17}
\p_\nu\sqrt{det(G)}(h_0)=\frac{tr(\nu)}{2}\quad,
\ee
and in the other hand
\be
\label{II.18}
\p^2_{\nu^2}\sqrt{det(G)}(h^0)=  e^{-2\la}\,det(\nu)-4^{-1}\ e^{-2\la}\,(tr(\nu))^2
\ee
Combining now (\ref{II.13}) and (\ref{II.17}) gives in one hand
\be
\label{II.19}
J_0^t\,\p_{\nu}J(h^0)=- e^{-2\la} \ \nu^0
\ee
where $\nu^0$ is the trace free part of $\nu$
\[
\nu^0:=\nu-\frac{tr(\nu)}{2}\, I_2\quad,
\]
and in the other hand
\be
\label{II.20}
J_0^t\,\p^2_{\nu^2}J(h^0)=2\ e^{-4\la}\ \nu^0\,\nu+e^{-4\la}\ \lf(det(\nu)-\lf(\frac{tr(\nu)}{2}\rg)^2\rg)\ I_2\quad.
\ee
Observe that
\be
\label{II.20zza}
det(\nu^0)=det(\nu)-\lf(\frac{tr(\nu)}{2}\rg)^2
\ee
This yields
\be
\label{II.21zza}
J_0^t\,\p^2_{\nu^2}J(h^0)=2\ e^{-4\la}\ \nu^0\,\nu+e^{-4\la}\ det(\nu^0)\ I_2\quad.
\ee
Inserting these expressions in (\ref{II.8a})...(\ref{II.9b}) gives
\be
\label{II.18a}
\p_\nu \Pi^0_{lk}(h^0)=-i\,\int_{D^2} e^{-2\la}\,(A^l)^t\, \nu^0\, X^k_0\ dx^2\quad,
\ee
and
\be
\label{II.19a}
\p_\nu \Pi^1_{lk}(h^0)=-i\,\int_{D^2}e^{-2\la}\,(B^l)^t\, \nu^0\, X^k_0\ dx^2\quad.
\ee
Using intrinsic notations with $\nu=\nu_{ij}\ dx_i\otimes dx_j$ and $\nu_0=\nu-2^{-1}\,tr_{h^0}\nu\ h^0$ and the fact that $<dx_i\otimes dx_j,dx_i\otimes dx_j>_{h^0}=e^{-4\la}$, we get
\be
\label{a.II.18a}
\p_\nu \Pi^0_{lk}(h^0)=-i\,\int_{\Sigma}<\eta_{a_l}\otimes\al^k, \nu-2^{-1}\,tr_{h^0}\nu\ h^0>_{h^0}\ dvol_{h^0}      \quad,
\ee
and
\be
\label{b.II.18a}
\p_\nu \Pi^1_{lk}(h^0)=-i\,\int_{\Sigma}<\eta_{b_l}\otimes\al^k, \nu-2^{-1}\,tr_{h^0}\nu\ h^0>_{h^0}\ dvol_{h^0}      \quad.
\ee
For the same reason the second derivatives of $\Pi$ at $h^0$ are given by
\be
\label{II.18b}
\begin{array}{l}
\ds\p^2_{\nu^2} \Pi^0_{lk}(h^0)=2i\,\int_{D^2}e^{-4\la}\ (A^l)^t\, \nu^0\,\nu\, X^k_0\ dx^2\\[5mm]
\ds +i\,\int_{D^2}  e^{-4\la}\ det(\nu^0)\ (A^l)^t\, X^k_0\ dx^2 -2i\,\int_{D^2} e^{-2\la}\ (A^l)^t\, \nu^0\, \nabla(\p_\nu \varphi^k(h^0))\ dx^2\quad,
\end{array}
\ee
and
\be
\label{II.19b}
\begin{array}{l}
\ds\p^2_{\nu^2} \Pi^1_{lk}(h^0)=2i\,\int_{D^2}e^{-4\la}\ (B^l)^t\, \nu^0\,\nu\, X^k_0\ dx^2\\[5mm]
\ds +i\,\int_{D^2}  e^{-4\la}\ det(\nu^0)\ (B^l)^t\, X^k_0\ dx^2 -2i\,\int_{D^2} e^{-2\la}\ (B^l)^t\, \nu^0\, \nabla(\p_\nu \varphi^k(h^0))\ dx^2\quad.
\end{array}
\ee

\subsection{The first and second derivatives of the Periods at an isothermic surface.}

\subsubsection{Weak isothermic immersions.}
We consider $\vec{\Phi}$ to be a weak, possibly branched, immersion in ${\mathcal F}_\Sigma$. As explained in the introduction 
such a weak, possibly branched, immersion defined a unique smooth Teichm\"uller class $\tau(\vec{\Phi}^\ast g_{{\R}^m})$
that we will simply denote $\tau_{\vec{\Phi}}$. In the computation below we will assume that $\vec{\Phi}$ is a global, possibly
branched, weak isothermic immersion (see \cite{Ri3} and \cite{Ri5}), that is to say there exists a non zero holomorphic quadratic differential\footnote{This is an holomorphic section of  $\wedge^{1-0}\Sigma\otimes\wedge^{1-0} \Sigma$.}
$q=f(z)\ dz\otimes dz$ such that
\be
\label{II.20a}
\Im(q,\vec{\mathfrak{h}}_0)_{wp}\equiv0\quad.
\ee
where $\vec{\mathfrak{h}}_0$ is the {\it Weingarten form} of $\vec{\Phi}$ defined as being the
section of ${\R}^m\otimes\wedge^{1-0}\Sigma\otimes\wedge^{1-0} \Sigma$  whose expression  in an arbitrary choice of complex coordinates is given by:
\be
\label{II.21}
\begin{array}{l}
\vec{\mathfrak{h}}_0:= 2\,\pi_{\vec{n}}(\p^2_{z^2}\vec{\Phi})\ dz\otimes dz\\[5mm]
\ds\quad=2^{-1}\,\pi_{\vec{n}}\lf(\p^2_{x^2_1}\vec{\Phi}-\p^2_{x^2_2}\vec{\Phi}-2\, i\ \p^2_{x_1x_2}\vec{\Phi}\rg)\ \ {dz\otimes dz}
\end{array}
\ee
where $\pi_{\vec{n}}$ is the orthogonal projection onto the plane orthogonal to $\vec{\Phi}_\ast T\Sigma$, $z=x_1+i x_2$ and $\p_z:={2^{-1}}\ [\p_{x_1}-i\p_{x_2}]$, moreover $(\cdot,\cdot)_{wp}$ is the following pointwise
hermitian product depending on the conformal immersion $\vec{\Phi}$
\[
(\psi_1\ dz\otimes dz,\psi_2\ dz\otimes dz)_{wp}:= e^{-4\la}\ \ov{\psi_1(z)}\ \psi_2(z)
\]
where $e^\la:=|\p_{x_1}\vec{\Phi}|=|\p_{x_2}\vec{\Phi}|$. For 2 sections $\Psi_1$, $\Psi_2$ of $\wedge^{0,1}\Sigma\otimes\wedge^{0,1}\Sigma$ we finally
denote
\be
\label{II.210az}
(\Psi_1,\Psi_2)_{WP}:=\int_{\Sigma}(\Psi_1,\Psi_2)_{wp}\ dvol_{h^0}
\ee
and $\Im(\Psi_1,\Psi_2)_{WP}$ realizes a scalar product on $\wedge^{0,1}\Sigma\otimes\wedge^{0,1}\Sigma$.

\medskip

Choosing complex coordinates in which $f(z)\equiv 1$ (this is possible away from the zeros of $q$) identity (\ref{II.20a}) becomes
(see \cite{Ri3} and \cite{Ri5})
\be
\label{II.21az}
\Im\lf(\vec{H}_0\rg)=2\Im\lf(\p_z\lf[e^{-2\la}\,\p_z\vec{\Phi}\rg]\rg)=0
\ee
where
\[
\vec{\mathfrak{h}}_0=e^{2\la}\,\vec{H}_0\ dz\otimes dz
\]
which is equivalent to the existence of $\vec{L}\in W{1,\infty}(D^2,{\R}^m)$ such that
\be
\label{II.21bz}
\lf\{
\begin{array}{l}
\p_{x_1}\vec{L}=e^{-2\la}\p_{x_1}\vec{\Phi}\quad\\[5mm]
\p_{x_2}\vec{L}=-e^{-2\la}\p_{x_2}\vec{\Phi}\quad
\end{array}
\rg.
\ee
If there would exists two \underbar{real} linearly independent holomorphic quadratic forms on $\Sigma$ such that (\ref{II.20a}) holds
then we would have away from isolated points
\[
\Re\lf(\vec{H}_0\rg)=2\,\Re\lf(\p_z\lf[e^{-2\la}\,\p_z\vec{\Phi}\rg]\rg)=0
\]
from which we would deduce
\[
\vec{\mathfrak{h}}_0\equiv 0\quad,
\]
which is equivalent for $(\Sigma,\vec{\Phi}^\ast g_{{\R}^m})$ to be {\it umbilic} that would contradicts the fact $genus(\Sigma)>0$.
Hence we have the following proposition
\begin{Prop}
\label{pr-II.2}
Let $\Sigma$ be a closed two-manifold with positive genus. Let $\vec{\Phi}$ be a global possibly branched weak isothermic immersion of ${\mathcal F}_\Sigma$. Assume we have two distinct holomorphic quadratic form $q_1$ and $q_2$ solving (\ref{II.20a})
then they are \underbar{real} linearly dependent : there exist a non trivial pair of real numbers $(t_1,t_2)\ne (0,0)$ s.t.
\[
t_1q_1+t_2 q_2=0\quad.
\]
\hfill $\Box$
\end{Prop}
In \cite{Ri3} we proved that the map ${\mathcal C}$ from ${\mathcal E}_\Sigma$, viewed as a Banach manifold over $W^{2,2}\cap W^{1,\infty}$,
into $T(\Sigma)$ which to $\vec{\Phi}$ assigns it's Teichm\"uller class $\tau(\vec{\Phi}^\ast g_{{\R}^m})$ is $C^1$ and it's differential is given by
\be
\label{II.22}
\p_{\vec{w}}{\mathcal C}({\Phi})=\sum_{j=1\cdots N}q_j\ \int_{\Sigma}\Im(q_j,\vec{h}_0)_{wp}\cdot\vec{w}\ dvol_{g_{\vec{\Phi}}}=\sum_{j=1\cdots N}q_j\ \Im(q_j,\vec{\mathfrak{h}}^0\cdot\vec{w})_{WP}
\ee
where $(q_j)_{j=1\cdots Q}$ is an orthonormal basis for the Weil-Peterson hermitian product on the space of Holomorphic quadratic form which identifies with the tangent space to the Teichm\"uller space $T(\Sigma)$ (i.e $Q=1$ for $g=1$ and
$Q=3g-3$ otherwise).

\medskip

On ${\mathcal F}_\Sigma$ we have that ${\mathcal C}$ is Frechet differentiable for perturbations $\vec{w}$ supported
on compact set which do not contain the branch points and the Frechet differential is also given by (\ref{II.22}).

\medskip

We now interpret the space of holomorphic quadratic differentials as being a \underbar{real} vector space
generated by $((q_j)_{j=1\cdots Q},(i\,q_j)_{j=1\cdots Q})$. Assuming $\vec{\Phi}$ is isothermic we have (\ref{II.20a})
or in other words there is $2Q-$tuple of reals $(t_1\cdots t_Q,s_1\cdots s_Q)\ne(0\cdots 0)$ such that
\be
\label{II.23}
\lf(\sum_{j=1}^Q t_j\ q_j+ s_j\ i\, q_j, \vec{h}_0\rg)_{WP}\equiv 0
\ee
This implies
\be
\label{II.23a}
\forall \vec{w}\in W^{2,2}\cap W^{1,\infty}\quad\quad\sum_{j=1}^Q t_j\ \int_{\Sigma}\Im(q_j,\vec{h}_0)_{wp}\cdot \vec{w}\ dvol_{g_{\vec{\Phi}}}+ \sum_{j=1}^Q s_j\ i\, \int_{\Sigma}\Im(q_j,\vec{h}_0)_{wp}\cdot \vec{w}\ dvol_{g_{\vec{\Phi}}}=0
\ee
or in other words $\p{\mathcal C}$  is included in the hyperplane of $T(\Sigma)$ with normal vector  $(t_1\cdots t_Q,s_1\cdots s_Q)$ in the orthonormal basis  $((q_j)_{j=1\cdots Q},(i\,q_j)_{j=1\cdots Q})$. If the range of $\p{\mathcal C}$ would be included in a strict subspace to this hyperplane there would be  another linearly independent family of ${\R}^{2Q}$, non parallel
to $(t_1\cdots t_Q,s_1\cdots s_Q)$ satisfying (\ref{II.23a}) and hence (\ref{II.23}) but this would contradict proposition~\ref{pr-II.2}.
Thus we have proved the following proposition.
\begin{Prop}
\label{pr-II.3}
Let $\Sigma$ be a closed two-manifold with positive genus. Let $\vec{\Phi}$ be a global possibly branched weak isothermic immersion of ${\mathcal F}_\Sigma$ then the real rank of the differential of the map which to $\vec{\Phi}$ assigns its Teichm\"uller class is exactly equal to $dim_{{\R}}T(\Sigma)-1$.\hfill $\Box$
\end{Prop}

\subsubsection{Computation of $\p_{\vec{w}}\Pi(\vec{\Phi})$ }

We consider a weak, possibly branched, conformal isothermic immersion $\vec{\Phi}$ and perturbations $\vec{w}\in W^{1,\infty}\cap W^{2,2}$ supported in a disc in $\Sigma$ that avoid the branched points of $\vec{\Phi}$, as well as the zeros of any non-trivial holomorphic quadratic form satisfying (\ref{II.20a}). Moreover we choose the complex coordinates in this disc in such a way
that $q= dz\otimes dz$ therefore (\ref{II.21az}) and (\ref{II.21bz}) hold.

\medskip

The perturbed metric obtained by adding $\vec{w}$ to $\vec{\Phi}$ is in these coordinates given by
\[
h(\vec{\Phi}+\vec{w})=h^0+\nu(\vec{\Phi}+\vec{w})=e^{2\la}\, I_2+(\nabla \vec{w})\,(\nabla \vec{\Phi})^t+(\nabla \vec{\Phi})\,(\nabla \vec{w})^t+(\nabla \vec{w})\,(\nabla \vec{w})^t\quad.
\]
Thus
\be
\label{II.24}
\lf\{
\begin{array}{l}
\ds\p_{\vec{w}}\nu(\vec{\Phi})=(\nabla \vec{w})\,(\nabla \vec{\Phi})^t+(\nabla \vec{\Phi})\,(\nabla \vec{w})^t\\[5mm]
\ds\p^2_{\vec{w}^2}\nu(0)=2(\nabla \vec{w})\,(\nabla \vec{w})^t
\end{array}
\rg.
\ee
and
\be
\label{II.25a}
\begin{array}{l}
\ds\p_{\vec{w}}\nu^0(\vec{\Phi})=\p_{\vec{w}}\nu(\vec{\Phi})-\frac{\p_{\vec{w}} tr(\nu)(\vec{\Phi})}{2}\,I_2\\[5mm]
\ds\quad=
\lf(
\begin{array}{cc}
\ds \p_{x_1}\vec{\Phi}\cdot\p_{x_1}\vec{w}-   \p_{x_2}\vec{\Phi}\cdot\p_{x_2}\vec{w}&\p_{x_1}\vec{\Phi}\cdot\p_{x_2}\vec{w}+ \p_{x_2}\vec{\Phi}\cdot\p_{x_1}\vec{w}   \\[5mm]
\ds  \p_{x_1}\vec{\Phi}\cdot\p_{x_2}\vec{w}+ \p_{x_2}\vec{\Phi}\cdot\p_{x_1}\vec{w}  &-\p_{x_1}\vec{\Phi}\cdot\p_{x_1}\vec{w}+ \p_{x_2}\vec{\Phi}\cdot\p_{x_2}\vec{w}
\end{array}
\rg)
\end{array}
\ee
We can assume \underbar{for the first derivative} that the perturbation $\vec{w}$ is normal to the surface. Then we have
\be
\label{II.25b}
\begin{array}{l}
\ds\sum_{i,j=1}^2\p_{\vec{w}}\nu^0(\vec{\Phi})_{ij}\ dx_i\otimes dx_j=\p_{\vec{w}}h-2^{-1}\ tr_{h^0}\p_{\vec{w}}h\\[5mm]
\ds=-\vec{w}\cdot\lf[\p^2_{x_1^2}\vec{\Phi}-\p^2_{x_2^2}\vec{\Phi}\rg]\ [dx_1^2-dx_2^2]-2\,\vec{w}\cdot\p^2_{x_1x_2}\vec{\Phi}\ [dx_1\otimes dx_2+dx_2\otimes dx_1]\\[5mm]
=-2\ \Re\lf(\vec{{\mathfrak{h}}}^0\cdot\vec{w}\rg)
\end{array}
\ee
Combining this identity with (\ref{a.II.18a}) and (\ref{b.II.18a})
\be
\label{a.II.18b}
\p_{\vec{w}} \Pi^0_{lk}(\vec{\Phi})=2\,i\,\int_{\Sigma}\vec{w}\cdot\Re<\eta_{a_l}\otimes\al^k,\vec{{\mathfrak{h}}}^0>_{\vec{\Phi}^\ast g_{{\R}^m}}\ dvol_{\vec{\Phi}^\ast g_{{\R}^m}}      \quad,
\ee
and
\be
\label{b.II.18b}
\p_{\vec{w}} \Pi^1_{lk}(\vec{\Phi})=2\,i\,\int_{\Sigma}\vec{w}\cdot\Re<\eta_{b_l}\otimes\al^k, \vec{{\mathfrak{h}}}^0>_{\vec{\Phi}^\ast g_{{\R}^m}}\ dvol_{\vec{\Phi}^\ast g_{{\R}^m}}    \quad.
\ee
where we are denoting $\Pi^0_{lk}(\vec{\Phi})$  or $\Pi^1_{lk}(\vec{\Phi})$ for $(\Pi^0_{lk}\circ h)(\vec{\Phi})$ and for resp. $(\Pi^1_{lk}\circ h)(\vec{\Phi})$.

\subsubsection{Free families of periods and isothermic surfaces.}

\begin{Dfi}
\label{df-free}
A sub-familly of periods $(\Pi^0_{kl},\Pi^1_{pq})_{(k,l)\in K\,,\, (p,q)\in P}$ where $K$ and $P$ are subsets of $\{1\cdots g\}^2$ is said to be ${\R}-$free
at a metric $h^0$ if the corresponding  differentials 
\[
(\p_\nu\Pi^0_{kl},\p_\nu\Pi^1_{pq})_{(k,l)\in K\,,\, (p,q)\in P}
\]
are independent for the ${\R}$-vector space structure- i.e. realizes a free family. The definition extends naturally to free families of linear combinations of periods. \hfill $\Box$
\end{Dfi}

Let $\eta_{a_l}$ be the harmonic 1 form representing the Poincar\'e dual of $a_l$ and let $\al^k$ be the family of harmonic 1-form on
$(\Sigma,h_0)$ satisfying (\ref{0.1}), (\ref{0.2}) and (\ref{0.3}). Recall the notations $\om^k=\al^k+i\ast_{0}\al^k$ for the corresponding basis
of holomorphic 1-forms and introduce $\sigma^l=i\,(\eta_{a_l}+i\,\ast_0\,\eta_{a_l})$ to be the holomorphic 1-forms associated to $\eta_{a_l}$.
We have already denoted $\al^k=X^k_1\,dx_1+X^k_2\, dx_2$ and $\eta_{a_l}=A^l_1\,dx_1+A^l_2\,dx_2$. Observe that we have
locally in conformal charts for $\ast_0$
\[
\om^k=(X^k_1-iX^k_2)\, dz\quad\quad\mbox{ and }\quad\quad\sigma^l=i\ (A^l_1-iA^l_2)\, dz\quad.
\]
Hence
\[
\om^k\otimes\sigma^l=i\,[(X^k_1\,A^l_1+X^k_2\,A^l_2)-i\,(X^k_1\,A^l_2+X^k_2\,A^l_1)]\ dz^2\quad.
\]
Let $\nu$ be a variation of the metric : $h=h^0+t\,\nu$. $\nu$ is an arbitrary symmetric tensor of $T^\ast\Sigma\otimes T^\ast\Sigma$. Denotes locally
in a given conformal charts $\nu:=\sum_{i,j=1}^2\nu_{ij}\ dx^i\otimes dx^j$. A short computation gives
\[
\nu^0=\nu-2^{-1}\,tr_{h^0}\nu\ h^0=e^{2\la}\ \Re(v^0\ dz^2)
\]
where $v^0=v^0_{\Re}+i\,v^0_{\Im}$ and $v^0_{\Re}=2^{-1}e^{-2\la}\,(\nu_{11}-\nu_{22})$ and $v^0_{\Im}=e^{-2\la}\nu_{12}$.
Denote also $\mu^0:=e^{2\la}\ v^0\ dz^2$. Observe that $\mu^0$ is an \underbar{arbitrary section} of the bundle of $1-0\otimes 1-0$ forms
as $\nu$ describes a neighborhood of metrics around $h^0$. Locally in these coordinates, we have in one hand
\be
\label{II.a25}
\begin{array}{l}
\ds\Im(\om^k\otimes\sigma^l,\mu^0)_{wp}=e^{-2\la}\ \Re\lf([(X^k_1\,A^l_1+X^k_2\,A^l_2)+i\,(X^k_1\,A^l_2+X^k_2\,A^l_1)]\ [v^0_{\Re}+i\,v^0_{\Im}]\rg)\\[5mm]
\quad= e^{-2\la}[v^0_{\Re}\,(X^k_1\,A^l_1+X^k_2\,A^l_2)-v^0_{\Im}\,(X^k_1\,A^l_2+X^k_2\,A^l_1)]
\end{array}
\ee
In the other hand we have
\be
\label{II.b25}
\begin{array}{l}
\ds<\eta_{a_l}\otimes\al^k,\nu^0>_{\vec{\Phi}^\ast g_{{\R}^m}}=\Re<\eta_{a_l}\otimes\al^k,\mu^0>_{\vec{\Phi}^\ast g_{{\R}^m}}\\[5mm]
\ds\quad=\sum_{i,j=1}^2\ e^{2\la}\ X^k_i\,A^l_j \lf[v^0_{\Re}\,\Re<dx_i\, dx_j, dz^2>-v^0_{\Im}\, \Im<dx_i\, dx_j, dz^2>\rg]\\[5mm]
\ds\quad=e^{-2\la}[v^0_{\Re}\,(X^k_1\,A^l_1+X^k_2\,A^l_2)-v^0_{\Im}\,(X^k_1\,A^l_2+X^k_2\,A^l_1)]
\end{array}
\ee
So we have proved that 
\be
\label{II.bb25}
\Im(\om^k\otimes\sigma^l,\mu^0)_{wp}=<\eta_{a_l}\otimes\al^k,\nu^0>_{h^0}\quad.
\ee
Similarly, denoting $\tau^l:=i\,(\eta_{b_l}+i\,\ast_0\,\eta_{b_l})$, we have
\be
\label{II.bc25}
\Im(\om^k\otimes\tau^l,\mu^0)_{wp}=<\eta_{b_l}\otimes\al^k,\nu^0>_{h^0}\quad.
\ee
Now, using (\ref{0.1}), we obtain that $$\eta_{a_l}=\ast\al^l$$ and then we have that $\sigma^l=\om^l$. Moreover using (\ref{0.0})...(\ref{II.2}) we have that
\[
\eta_{b_l}=\al_l+\sum_{j=1}^gc_{lj}\ \ast_0\,\al_j
\]
which gives
\[
\tau^l=i\,(\al^l+i\,\ast_0\,\al^l)+\sum_{j=1}^gc_{lj}\ (\al^l+i\,\ast_0\,\al^l)=i\,\om^l+\sum_{j=1}^gc_{lj}\ \om^j\quad.
\]
Hence we deduce
\be
\label{II.bd25}
<\eta_{a_l}\otimes\al^k,\nu^0>_{h^0}=\Im(\om^k\otimes\om^l,\mu^0)_{wp}
\ee
and
\be
\label{II.be25}
<\eta_{b_l}\otimes\al^k,\nu^0>_{h^0}=\Re(\om^k\otimes\om^l,\mu^0)_{WP}+\sum_{j=1}^gc_{lj}\,\Im(\om^k\otimes\om^j,\mu^0)_{WP}
\ee
Combining (\ref{II.bd25}) and (\ref{II.be25}) with (\ref{a.II.18a}) and (\ref{b.II.18a})

\be
\label{II.be26}
\p_\nu \Pi^0_{lk}(h^0)=-i\,\Im(\om^k\otimes\om^l,\mu^0)_{WP}     \quad,
\ee
and
\be
\label{II.be27}
\p_\nu \Pi^1_{lk}(h^0)=-i\,\Im\lf(i\, \om^k\otimes\om^l+\sum_{j=1}^gc_{jl}\ \om^k\otimes\om^j,\mu^0\rg)_{WP}     \quad.
\ee
Again, since by varying $\nu$, $\mu^0$ describes the full space of sections of $(1-0)^2-$forms we obtain the following proposition

\begin{Prop}
\label{pr-free}
A sub-familly of periods $(\Pi^0_{kl},\Pi^1_{pq})_{(k,l)\in K\,,\, (p,q)\in P}$ where $K$ and $P$ are subsets of $\{1\cdots g\}^2$ is  ${\R}-$free
if and only if the following family of holomorphic quadratic forms are ${\R}-$independent
\[
q^0_{kl}:=\om^k\otimes\om^l \quad \mbox{where }{(k,l)\in K}\quad\mbox{together with }\quad q^1_{pq}:=i\, \om^p\otimes\om^q+\sum_{j=1}^gc_{jq}\ \om^p\otimes\om^j\quad \mbox{where }{(p,q)\in P}
\]
\hfill $\Box$
\end{Prop}
Applying (\ref{II.bd25}) and (\ref{II.be25}) to $\mu^0=\mathfrak{h}^0_i$ for $i=1\cdots m$ we have using (\ref{a.II.18b}) and (\ref{b.II.18b})
\be
\label{II.e25}
\p_{\vec{w}} \Pi^0_{lk}(\vec{\Phi})=2\,i\,\int_{\Sigma}\vec{w}\cdot\Im(\om^k\otimes\om^l,\vec{\mathfrak{h}}_0)_{wp}\ dvol_{\vec{\Phi}^\ast g_{{\R}^m}} 
=2\,i\, (q^0_{kl}, \vec{\mathfrak{h}}_0)_{WP}    \quad,
\ee
and
\be
\label{II.f25}
\begin{array}{l}
\ds\p_{\vec{w}} \Pi^1_{lk}(\vec{\Phi})=2\,i\,\int_{\Sigma}\vec{w}\cdot\Im(i\ \om^k\otimes\om^l,\vec{\mathfrak{h}}_0)_{wp}\ dvol_{\vec{\Phi}^\ast g_{{\R}^m}} \\[5mm]
\ds\quad\quad+2\,i\,\sum_{j=1}^g c_{jl} \int_{\Sigma}\vec{w}\cdot\Im(\om^k\otimes\om^j,\vec{\mathfrak{h}}_0)_{wp}\ dvol_{\vec{\Phi}^\ast g_{{\R}^m}}=2\,i\, (q^1_{kl}, \vec{\mathfrak{h}}_0)_{WP} \quad.
\end{array}
\ee
We Deduce then from these expressions the following proposition
\begin{Prop}
\label{pr-free2}
Let $\Phi$ be a weak immersion from ${\mathcal F}_\Sigma$.
Let $(\Pi^0_{kl},\Pi^1_{pq})_{(k,l)\in K\,,\, (p,q)\in P}$ where $K$ and $P$ are subsets of $\{1\cdots g\}^2$ be a  ${\R}-$free sub-familly of periods
at $h^0=\vec{\Phi}^\ast g_{{\R}^m}$. Assume
 \[
(\p(\Pi^0_{kl}\circ h)(\vec{\Phi}),\p(\Pi^1_{pq}\circ h)(\Phi))_{(k,l)\in K\,,\, (p,q)\in P}\quad\mbox{ is not free }\quad
\]
as pure imaginary one forms on $W^{2,2}\cap W^{1,\infty}(\Sigma,{\R}^m)$ then
\be
\label{II.g25}
rk\lf(\p(\Pi^0_{kl}\circ h)(\vec{\Phi}),\p(\Pi^1_{pq}\circ h)(\Phi))_{(k,l)\in K\,,\, (p,q)\in P}\rg)=|K|+|P|-1
\ee
and $\vec{\Phi}$ is isothermic. This proposition extends naturally to a free family of linear combinations of periods.\hfill $\Box$
\end{Prop}
\noindent{\bf Proof of proposition~\ref{pr-free2}.}
Under the assumption of the proposition, because of (\ref{II.e25}) and (\ref{II.f25}), there exists a non zero holomorphic quadratic differential $q$ such that
\be
\label{II.g26}
\Im(q,\vec{\mathfrak{h}}^0)\equiv 0\quad,
\ee
which is equivalent to the fact that $\vec{\Phi}$ is isothermic. Because of proposition~\ref{pr-II.2} the number of independent holomorphic quadratic differential 
such that (\ref{II.g26}) is at most 1 hence the decrease of the dimension between 
$\lf(\p\Pi^0_{kl}(h^0),\p\Pi^1_{pq}(h^0)\rg)_{(k,l)\in K\,,\, (p,q)\in P}$ and
$\lf(\p(\Pi^0_{kl}\circ h)(\vec{\Phi}),\p(\Pi^1_{pq}\circ h)(\Phi)\rg)_{(k,l)\in K\,,\, (p,q)\in P}$ is at most 1. This concludes the proof of proposition~\ref{pr-free2}.
\hfill $\Box$

\medskip

A classical theorem of M.Noether asserts that if the genus of $\Sigma$ satisfies $g\le2$ or if the conformal class on $\Sigma$ defined
by
$\vec{\Phi}^\ast g_{{\R}^m}$ is non hyper-elliptic\footnote{Recall that a conformal structure $(\Sigma,c)$ is hyper-elliptic if the subspace
of tensor products of holomorphic 1-forms in the space of holomorphic quadratic differentials is of real dimension $4g-2$ and, for $g>2$ this subspace is a complex analytic submanifold of dimension $2g-1$ to the Teichm\"uller space (see for instance \cite{Na} theorem 4.1.4). }  then the family $\om^k\otimes\om^l$ is of real dimension $6g-6$ and generates the space of holomorphic quadratic differentials. Then the real space generated by $\sum_{l,k}d\Pi_{l,k}\ \om^k\otimes\om^l$ is equal to the real space underlying the complex space generated by $\p{\mathcal C}$ in the tangent
space to the Teichm\"uller space at $\vec{\Phi}^\ast g_{{\R}^m}$. Hence we have the following partial reciproque of proposition~\ref{pr-free}
\begin{Prop}
\label{pr-isot}
Let $\vec{\Phi}$ be a weak immersion from ${\mathcal F}_\Sigma$. Assume either $g\le 2$ or assume that $(\Sigma, \vec{\Phi})$ is not hyper-elliptic if $g>2$.
Then $\vec{\Phi}$ is isothermic if and only if there exists a free family of periods $(\Pi^0_{kl},\Pi^1_{pq})_{(k,l)\in K\,,\, (p,q)\in P}$ such that
$(\p(\Pi^0_{kl}\circ h)(\vec{\Phi}),\p(\Pi^1_{pq}\circ h)(\Phi))_{(k,l)\in K\,,\, (p,q)\in P}$ is not free.\hfill $\Box$
\end{Prop}

\subsubsection{The computation of the second derivative of the period matrix at an isothermic surface.}

We are now assuming that $\vec{\Phi}$ is isothermic and we are working in local conformal coordinates such that (\ref{II.21bz}) holds. Rewriting 
(\ref{II.25a}) in terms of $\vec{L}$ gives

\be
\label{II.25}
\ds\p_{\vec{w}}\nu^0(\vec{\Phi})=
e^{2\la}\lf(
\begin{array}{cc}
\ds \nabla\vec{L}\cdot\nabla\vec{w}&\nabla\vec{w}\cdot\nabla^\perp\vec{L}  \\[5mm]
\ds  \nabla\vec{w}\cdot\nabla^\perp\vec{L} & -\nabla\vec{L}\cdot\nabla\vec{w}     \end{array}
\rg)= e^{2\la}\ \nabla \vec{w}\cdot\nabla\vec{L}\ S_1 +e^{2\la}\ \nabla \vec{w}\cdot\nabla^\perp\vec{L}\ S_2
\ee
where
\[
S_1:=\lf(
\begin{array}{cc}
1&0\\[5mm]
0&-1
\end{array}\rg)\quad\mbox{ and }\quad
S_2:=\lf(
\begin{array}{cc}
0&1\\[5mm]
1&0
\end{array}\rg)
\]
This implies that
\be
\label{II.26}
\p_{\vec{w}}\lf(\Pi^0_{lk}\circ h\rg)(\vec{\Phi})=-i\,\int_{D^2}\nabla\vec{L}\cdot\nabla\vec{w}\  (A^l)^t\,S_1\, X^k_0\ dx^2-i\,\int_{D^2}\nabla^\perp\vec{L}\cdot\nabla\vec{w}\   (A^l)^t\,S_2\, X^k_0\ dx^2\quad,
\ee
and
\be
\label{II.27}
\p_{\vec{w}}\lf(\Pi^1_{lk}\circ h\rg)(\vec{\Phi})=-i\,\int_{D^2}\nabla\vec{L}\cdot\nabla\vec{w}\  (B^l)^t\,S_1\, X^k_0\ dx^2-i\,\int_{D^2}\nabla^\perp\vec{L}\cdot\nabla\vec{w}\   (B^l)^t\,S_2\, X^k_0\ dx^2\quad,
\ee
We have
\be
\label{II.28}
\p^2_{\vec{w}^2}\lf(\Pi^0_{lk}\circ h\rg)(\vec{\Phi})=\p^2_{(\p_{\vec{w}}\nu)^2} \Pi^0_{lk}(h^0)+\p_{\p^2_{\vec{w}^2}\nu}\Pi^0_{lk}(h^0)
\ee
Combining (\ref{II.18a}), (\ref{II.18b}), (\ref{II.24}), (\ref{II.25}) and (\ref{II.28}) gives
\be
\label{II.29}
\begin{array}{l}
\ds\p^2_{\vec{w}^2}\lf(\Pi^0_{lk}\circ h\rg)(\vec{\Phi})=2i\,\int_{D^2}e^{-4\la}\ (A^l)^t\, \p_{\vec{w}}\nu^0\,\p_{\vec{w}}\nu\, X^k_0\ dx^2\\[5mm]
\ds +i\,\int_{D^2}  e^{-4\la}\ det(\p_{\vec{w}}\nu^0)\ (A^l)^t\, X^k_0\ dx^2 -2i\,\int_{D^2} e^{-2\la}\ (A^l)^t\, \nu^0\, \nabla(\p_{\p_{\vec{w}}} \varphi^k(h^0))\ dx^2\\[5mm]
\ds -i\,\int_{D^2} e^{-2\la}\,(A^l)^t\, \p^2_{{\vec{w}}^2}\nu^0\, X^k_0\ dx^2
\end{array}
\ee
We have
\be
\label{II.30}
(A^l)^t\, \p_{\vec{w}}\nu^0\,\p_{\vec{w}}\nu\, X^k_0=(A^l)^t\, \p_{\vec{w}}\nu^0\,\p_{\vec{w}}\nu^0\, X^k_0+\frac{tr(\p_{\vec{w}}\nu)}{2}\,(A^l)^t\, \p_{\vec{w}}\nu^0\, X^k_0\quad,
\ee
and
\be
\label{II.31}
\frac{tr(\p_{\vec{w}}\nu)}{2}=\nabla\vec{\Phi}\cdot\nabla\vec{w}\quad.
\ee
Thus
\be
\label{II.32}
\begin{array}{l}
\ds\int_{D^2}e^{-4\la}\ (A^l)^t\, \p_{\vec{w}}\nu^0\,\p_{\vec{w}}\nu\, X^k_0\ dx^2=\int_{D^2}e^{-4\la}\ (A^l)^t\, \p_{\vec{w}}\nu^0\,\p_{\vec{w}}\nu^0\, X^k_0\ dx^2\\[5mm]
\ds+\int_{D^2} \ov{\nabla}\vec{L}\cdot\nabla\vec{w}\ \nabla\vec{L}\cdot\nabla\vec{w}\ (A^l)^t\, S_1\, X^k_0\ dx^2+\int_{D^2} \ov{\nabla}\vec{L}\cdot\nabla\vec{w}\ \nabla^\perp\vec{L}\cdot\nabla\vec{w}\ (A^l)^t\, S_2\, X^k_0\ dx^2\quad,
\end{array}
\ee
where we denote
\[
\ov{\nabla}f=\lf(
\begin{array}{c} 
\p_{x_1}f\\[5mm]
-\p_{x_2}f
\end{array}
\rg)
\]
We have also
\be
\label{II.32a}
e^{-4\la}\ \p_{\vec{w}}\nu^0\,\p_{\vec{w}}\nu^0=\lf[((\nabla\vec{w}\cdot\nabla\vec{L})^2+(\nabla\vec{w}\cdot\nabla^\perp\vec{L})^2\rg]\ I_2\quad.
\ee
We have moreover using (\ref{II.25})
\be
\label{II.33}
e^{-4\la}\, det(\p_{\vec{w}}\nu^0)=-(\nabla\vec{w}\cdot\nabla\vec{L})^2-(\nabla\vec{w}\cdot\nabla^\perp\vec{L})^2
\ee
Hence, combining (\ref{II.24}), (\ref{II.29})...(\ref{II.33}) we obtain
\be
\label{II.34}
\begin{array}{l}
\ds\p^2_{\vec{w}^2}\lf(\Pi^0_{lk}\circ h\rg)(\vec{\Phi})=i\,\int_{D^2}\  \lf[((\nabla\vec{w}\cdot\nabla\vec{L})^2+(\nabla\vec{w}\cdot\nabla^\perp\vec{L})^2\rg]\    (A^l)^t\, X^k_0\ dx^2\\[5mm]
\ds+2i\int_{D^2} \ov{\nabla}\vec{L}\cdot\nabla\vec{w}\ \nabla\vec{L}\cdot\nabla\vec{w}\ (A^l)^t\, S_1\, X^k_0\ dx^2+2i\int_{D^2} \ov{\nabla}\vec{L}\cdot\nabla\vec{w}\ \nabla^\perp\vec{L}\cdot\nabla\vec{w}\ (A^l)^t\, S_2\, X^k_0\ dx^2\\[5mm]
\ds-2i\,\int_{D^2}\ (A^l)^t\, [\nabla \vec{w}\cdot\nabla\vec{L}\ S_1 +\ \nabla \vec{w}\cdot\nabla^\perp\vec{L}\ S_2]\, \nabla(\p_{\p_{\vec{w}}\nu} \varphi^k(h^0))\ dx^2\\[5mm]
\ds-i\,\int_{D^2} e^{-2\la}\,(A^l)^t\,\p^2_{\vec{w}^2}\nu^0 \, X^k_0\ dx^2
\end{array}
\ee
Denoting $a^k_{\vec{w}}:=\p_{\p_{\vec{w}}\nu} \varphi^k(h^0)-\ov{\p_{\p_{\vec{w}}\nu} \varphi^k(h^0)}$ where $\ov{\p_{\p_{\vec{w}}\nu} \varphi^k(h^0)}$ is the average
of $\p_{\p_{\vec{w}}\nu} \varphi^k(h^0)$ on $\Sigma$ for $h^0$, we have
\be
\label{II.34a}
\begin{array}{l}
\ds-2i\,\int_{D^2}\ (A^l)^t\, [\nabla \vec{w}\cdot\nabla\vec{L}\ S_1 +\ \nabla \vec{w}\cdot\nabla^\perp\vec{L}\ S_2]\, \nabla(\p_{\p_{\vec{w}}\nu} \varphi^k(h^0))\ dx^2\\[5mm]
\ds=-2i\int_{D^2} (A^l)^t\,V\, \p_{x_1}a^k_{\vec{w}}\ dx^2+2i\int_{D^2} (A^l)^t\,V^\perp\, \p_{x_2}a^k_{\vec{w}}\ dx^2\\[5mm]
\ds=2i\int_{D^2} a^k_{\vec{w}}\ \lf[\p_{x_1}(A^l)^t\ V-\p_{x_2}(A^l)^t\ V^\perp\rg] dx^2
+2i\int_{D^2} a^k_{\vec{w}}\ (A^l)^t\lf[\p_{x_1}V-\p_{x_2}V^\perp\rg]\ dx^2
\end{array}
\ee 
where
\[
V:=\lf(
\begin{array}{c}
\nabla\vec{w}\cdot\nabla\vec{L}     \\[5mm]
\nabla\vec{w}\cdot\nabla^\perp\vec{L}
\end{array}
\rg)
\quad\mbox{ and }\quad V^\perp:=\lf(
\begin{array}{c}
-\nabla\vec{w}\cdot\nabla^\perp\vec{L}     \\[5mm]
\nabla\vec{w}\cdot\nabla\vec{L}
\end{array}
\rg)\quad.
\]
Thus we have
\be
\label{II.34b}
\begin{array}{l}
\ds\p^2_{\vec{w}^2}\lf(\Pi^0_{lk}\circ h\rg)(\vec{\Phi})=i\,\int_{D^2}\  \lf[((\nabla\vec{w}\cdot\nabla\vec{L})^2+(\nabla\vec{w}\cdot\nabla^\perp\vec{L})^2\rg]\    (A^l)^t\, X^k_0\ dx^2\\[5mm]
\ds+2i\int_{D^2} \ov{\nabla}\vec{L}\cdot\nabla\vec{w}\ \nabla\vec{L}\cdot\nabla\vec{w}\ (A^l)^t\, S_1\, X^k_0\ dx^2+2i\int_{D^2} \ov{\nabla}\vec{L}\cdot\nabla\vec{w}\ \nabla^\perp\vec{L}\cdot\nabla\vec{w}\ (A^l)^t\, S_2\, X^k_0\ dx^2\\[5mm]
\ds+2i\int_{D^2} a^k_{\vec{w}}\ \lf[\p_{x_1}(A^l)^t\ V-\p_{x_2}(A^l)^t\ V^\perp\rg] dx^2
+2i\int_{D^2} a^k_{\vec{w}}\ (A^l)^t\lf[\p_{x_1}V-\p_{x_2}V^\perp\rg]\ dx^2\\[5mm]
\ds-i\,\int_{D^2} e^{-2\la}\,(A^l)^t\,\p^2_{\vec{w}^2}\nu^0\, X^k_0\ dx^2
\end{array}
\ee
From (\ref{II.1a})  we have, since $\varphi^k(h^0)=0$
\be
\label{II.35}
\Delta_h(\varphi^k(h))=\ast_h d(\ast_h\al_0^k)=\ast_h d\lf(\sum_{i,j=1}^2 J^{i,j}(h) X_j^k(h^0)\ dx_i\rg)
\ee
Since $\varphi^k(h^0)=0$ and since $d(\ast_{h^0}\al_0^k)=0$ we have\footnote{Recall that $\Delta_{h^0}$ is the positive
Laplace Beltrami operator for the metric $h^0$ given in the local complex charts by
\[
\Delta_{h^0}f=-e^{-2\la}\ \Delta f=-e^{-2\la}\ [\p^2_{x_1^2}f+\p^2_{x_2^2}f]\quad.
\]}
\be
\label{II.36}
-\Delta(\p_{\nu}\varphi^k(h^0))=\p_{x_1}\lf(\sum_{j=1}^2 \p_{\nu}J^{2,j}(h^0) X_j^k(h^0)\rg)-\p_{x_2}\lf(\sum_{j=1}^2 \p_\nu J^{1,j}(h^0) X_j^k(h^0)\rg)
\ee
Hence for any $f\in C^\infty(\Sigma)$ one has\footnote{Observe that the integrand $\nabla f\cdot\nabla a\ dx^2$ is independant
of the complex coordinate we locally choose.}
\be
\label{II.36a}
\int_{\Sigma}\nabla f\cdot\nabla a\ dx^2=\int_{D^2}\ (\nabla f)^t\ J_0\p_{\nu}J(h^0)\, X_0^k\ dx^2=
\int_{D^2}\ e^{-2\la}\, (\nabla f)^t\ \nu^0\, X_0^k\ dx^2
\ee
where we have used (\ref{II.19}). Using now (\ref{II.25}), we obtain
\be
\label{II.36b}
\int_{\Sigma}\nabla f\cdot\nabla a^k_{\vec{w}}\ dx^2=\int_{D^2}\lf[\p_{x_1}f\, V-\p_{x_2}f\, V^\perp\rg]^t\, X_0^k\ dx^2
\ee
or in other words
\be
\label{II.37a}
\Delta a^k_{\vec{w}}=\p_{x_1}(V^t\, X^k_0)-\p_{x_2}((V^\perp)^t\, X^k_0)\quad.
\ee
Let $b^k_{\vec{w}}$ be the function of average $0$ on $\Sigma$ for $h^0$ given by
\be
\label{II.37ab}
\ds\Delta b^k_{\vec{w}}=V^t\,\p_{x_1} X_0^k-(V^\perp)^t\,\p_{x_2} X_0^k
\ee
We have then
\be
\label{II.37ac}
\begin{array}{l}
\ds\Delta a^k_{\vec{w}}=\Delta b^k_{\vec{w}}+\lf[\p_{x_1}V-\p_{x_2}V^\perp\rg]^t\ X_0^k\\[5mm]
\quad=\Delta \lf(b^k_{\vec{w}}+\Delta^{-1}\lf[\p_{x_1}V-\p_{x_2}V^\perp\rg]^t\ X_0^k\rg)-
2\nabla\lf(\Delta^{-1}\lf[\p_{x_1}V-\p_{x_2}V^\perp\rg]^t\rg)\cdot\nabla X_k^0    \quad.
\end{array}
\ee
where we have used the fact that $\Delta X_k^0=0$. Denote
\[
c_{\vec{w}}^k:=b_{\vec{w}}^k-2\Delta^{-1}\lf(2\nabla\lf(\Delta^{-1}\lf[\p_{x_1}V-\p_{x_2}V^\perp\rg]^t\rg)\cdot\nabla X_k^0  \rg)\quad.
\]
so that we have
\be
\label{II.37aca}
a_{\vec{w}}^k=c_{\vec{w}}^k+\Delta^{-1}\lf[\p_{x_1}V-\p_{x_2}V^\perp\rg]^t\ X_0^k\quad.
\ee
We have
\be
\label{II.34d}
\begin{array}{l}
\ds\p^2_{\vec{w}^2}\lf(\Pi^0_{lk}\circ h\rg)(\vec{\Phi})=i\,\int_{D^2}\  \lf[((\nabla\vec{w}\cdot\nabla\vec{L})^2+(\nabla\vec{w}\cdot\nabla^\perp\vec{L})^2\rg]\    (A^l)^t\, X^k_0\ dx^2\\[5mm]
\ds+2i\int_{D^2} \ov{\nabla}\vec{L}\cdot\nabla\vec{w}\ \nabla\vec{L}\cdot\nabla\vec{w}\ (A^l)^t\, S_1\, X^k_0\ dx^2+2i\int_{D^2} \ov{\nabla}\vec{L}\cdot\nabla\vec{w}\ \nabla^\perp\vec{L}\cdot\nabla\vec{w}\ (A^l)^t\, S_2\, X^k_0\ dx^2\\[5mm]
\ds+2i\int_{D^2} (A^l)^t\lf[\p_{x_1}V-\p_{x_2}V^\perp\rg]\ \Delta^{-1}\lf[\p_{x_1}V-\p_{x_2}V^\perp\rg]^t\ X_0^k\ dx^2\\[5mm]
\ds+2i\int_{D^2} a^k_{\vec{w}}\ \lf[\p_{x_1}(A^l)^t\ V-\p_{x_2}(A^l)^t\ V^\perp\rg] \ dx^2
+2i\int_{D^2} c^k_{\vec{w}}\ (A^l)^t\lf[\p_{x_1}V-\p_{x_2}V^\perp\rg]\ dx^2\\[5mm]
\ds-i\,\int_{D^2} e^{-2\la}\,(A^l)^t\,\p^2_{\vec{w}^2}\nu^0 \, X^k_0\ dx^2
\end{array}
\ee
We write $Y^t=(y_1,y_2)=[\p_{x_1}V-\p_{x_2}V^\perp]=(\p_{x_1}v_1+\p_{x_2}v_2,\p_{x_1}v_2-\p_{x_2}v_1)$. Observe that
\be
\label{II.34e}
\begin{array}{l}
\ds 2\, Y\ \Delta^{-1}Y^t=(y_1\,\Delta^{-1}y_1+y_2\,\Delta^{-1}y_2)\ I_2+(y_1\,\Delta^{-1}y_1-y_2\,\Delta^{-1}y_2)\ S_1\\[5mm]
\ds +(y_1\,\Delta^{-1}y_2+y_2\,\Delta^{-1}y_1)\ S_2
+(-y_1\,\Delta^{-1}y_2+y_2\,\Delta^{-1}y_1)\ J_0
\end{array}
\ee
Denoting $y=y_1+i y_2$ and $v=v_1+iv_2$ we have $y=2\,\p_z v$, we have moreover
\be
\label{II.34f}
 2\, Y\ \Delta^{-1}Y^t=\Re(y\,\Delta^{-1}\ov{y})\, I_2+\Im(y\,\Delta^{-1}\ov{y})\, J_0+\Re(y\,\Delta^{-1}{y})\, S_1+\Im(y\,\Delta^{-1}{y})\, S_2
\ee
and hence
\be
\label{II.34g}
\begin{array}{l}
\ds 2 \lf[\p_{x_1}V-\p_{x_2}V^\perp\rg]\ \Delta^{-1}\lf[\p_{x_1}V-\p_{x_2}V^\perp\rg]^t= -|v|^2\, I_2-4\Re\lf( v\,\Delta^{-1}\p^2_{z^2}v\rg)\, S_1-4\Im\lf( v\,\Delta^{-1}\p^2_{z^2}v\rg)\, S_2\\[5mm]
\quad\quad\quad+4\,\Re\lf(\p_{z}\lf(v\,\Delta^{-1}\p_{\ov{z}}\ov{v}\rg)\rg)I_2+4\,\Im\lf(\p_{z}\lf(v\,\Delta^{-1}\p_{\ov{z}}\ov{v}\rg)\rg)J_0\\[5mm]
\quad\quad\quad+4\,\Re\lf( \p_{z}\lf(v\,\Delta^{-1}\p_{z}v\rg)\rg) S_1+4\,\Im\lf( \p_{z}\lf(v\,\Delta^{-1}\p_{z}v\rg)\rg)S_2
\end{array}
\ee
Observe that, because of (\ref{II.24}), we have
\be
\label{II.34gg}
\p^2_{\vec{w}^2}\nu^0=\lf(
\begin{array}{cc}
\nabla\vec{w}\cdot\ov{\nabla}\vec{w} & 2\,\p_{x_1}\vec{w}\cdot\p_{x_2}\vec{w} \\[5mm]
 2\,\p_{x_1}\vec{w}\cdot\p_{x_2}\vec{w} &-\nabla\vec{w}\cdot\ov{\nabla}\vec{w}
\end{array}
\rg)=\nabla\vec{w}\cdot\ov{\nabla}\vec{w}\, S_1+2\,\p_{x_1}\vec{w}\cdot\p_{x_2}\vec{w}\, S_2\quad.
\ee
Combining (\ref{II.34d}), (\ref{II.34g}) and (\ref{II.34gg}) gives
\be
\label{II.34h}
\begin{array}{l}
\ds\p^2_{\vec{w}^2}\lf(\Pi^0_{lk}\circ h\rg)(\vec{\Phi})=\\[5mm]
\ds+2i\int_{D^2} \lf[\ov{\nabla}\vec{L}\cdot\nabla\vec{w}\ \nabla\vec{L}\cdot\nabla\vec{w}-2\,\Re\lf( v\,\Delta^{-1}\p^2_{z^2}v\rg)
-4^{-1}\ |\nabla \vec{L}|^2\ \nabla\vec{w}\cdot\ov{\nabla} \vec{w} \rg](A^l)^t\, S_1\, X^k_0\ dx^2\\[5mm]
\ds+2i\int_{D^2} \lf[\ov{\nabla}\vec{L}\cdot\nabla\vec{w}\ \nabla^\perp\vec{L}\cdot\nabla\vec{w}-2\Im\lf( v\,\Delta^{-1}\p^2_{z^2}v\rg)- 2^{-1}|\nabla\vec{L}|^2\,\p_{x_1}\vec{w}\cdot\p_{x_2}\vec{w}\rg]\ (A^l)^t\, S_2\, X^k_0\ dx^2\\[5mm]
\ds+4i\,\int_{D^2}\Re\lf(\p_{z}\lf(v\,\Delta^{-1}\p_{\ov{z}}\ov{v}\rg)\rg)\,(A^l)^t\,X^k_0+4i\,\int_{D^2}\Im\lf(\p_{z}\lf(v\,\Delta^{-1}\p_{\ov{z}}\ov{v}\rg)\rg)\,(A^l)^t\,J_0X^k_0\ dx^2\\[5mm]
\ds+4i\,\int_{D^2}\Re\lf(\p_{z}\lf(v\,\Delta^{-1}\p_{{z}}{v}\rg)\rg)\,(A^l)^t\,S_1X^k_0+4i\,\int_{D^2}\Im\lf(\p_{z}\lf(v\,\Delta^{-1}\p_{{z}}{v}\rg)\rg)\,(A^l)^t\,S_2X^k_0\ dx^2\\[5mm]
\ds+2i\int_{D^2} a^k_{\vec{w}}\ \lf[\p_{x_1}(A^l)^t\ V-\p_{x_2}(A^l)^t\ V^\perp\rg] \ dx^2
+2i\int_{D^2} c^k_{\vec{w}}\ (A^l)^t\lf[\p_{x_1}V-\p_{x_2}V^\perp\rg]\ dx^2
\end{array}
\ee
where we have used the fact that $|\nabla\vec{L}|^2=2\, e^{-2\la}$ and we recall that 
\[
v=\nabla\vec{w}\cdot(\nabla\vec{L}+i\nabla^\perp\vec{L})=4\ \p_{\ov{z}}\vec{w}\cdot\p_z\vec{L}\quad.
\]
\subsubsection{Some special directions in the closure of the range of $\p^2(\Pi\circ h)$ for isothermic immersions.}  

The goal of this subsection is to establish the following theorem that will be deduced from the computations in the previous sections.

\begin{Prop}
\label{pr-second}
Let $\vec{\Phi}$ be an isothermic immersion of the surface $\Sigma$. Under the above notations, for any real function $\phi$ on ${\C}$,
and for almost every $x^0$ in $\Sigma$ there exists a sequence of directions $\vec{w}_\ep$ such that $supp\vec{w}_\ep\subset B_\ep(x^0)$,
such that $\|\nabla\vec{w}_\ep\|_{L^2}$ is uniformly bounded and for all $(l,k)\in \{1\cdots g\}^2$
\be
\label{II.450}
\p_{\vec{w}_\ep}\lf(\Pi^0_{lk}\circ h\rg)(\vec{\Phi})=O(\ep)\quad\mbox{ and }\quad\p_{\vec{w}_\ep}\lf(\Pi^1_{lk}\circ h\rg)(\vec{\Phi})=O(\ep)
\ee
and such that
\be
\label{II.520}
\lim_{\ep\rightarrow 0}\p^2_{\vec{w}^2_\ep}\lf(\Pi^0_{lk}\circ h\rg)(\vec{\Phi})=-4\,i\, \Im\lf(\om^l\otimes\om^k(x^0),\int_{{\C}}\p\phi\otimes\p\phi\rg)_{wp}\quad.
\ee
and 
\be
\label{II.530}
\lim_{\ep\rightarrow 0}\p^2_{\vec{w}^2_\ep}\lf(\Pi^1_{lk}\circ h\rg)(\vec{\Phi})=-4\,i \,\Im\lf(i\,\om^l\otimes\om^k(x^0)+\sum_{j=1}^gc_{kj}\,\om^l\otimes\om^j(x^0),\int_{{\C}}\p\phi\otimes\p\phi\rg)_{wp}\quad.
\ee

\hfill $\Box$
\end{Prop}

\noindent{\bf Proof of proposition~\ref{pr-second}.}

In the coordinates chart let
$$\phi_\ep(x)=\phi\lf(\frac{x-x^0}{\ep}\rg)$$ where the point $x^0=(x^0_1,x^0_2)$ is a Lebesgue point for $\nabla\vec{L}$ and $\phi(x)$ is a smooth function compactly supported in the unit ball $B_1(0)$. Let $\vec{\nu}$ a constant unit vector perpendicular to the surface at $x^0$. Since $x^0$ is chosen to be a Lebesgue point for $\nabla \vec{L}$ - which is tangent to the surface almost everywhere - we have
\be
\label{II.43}
\lim_{\ep\rightarrow 0}\frac{1}{\ep^2}\int_{D^2_\ep}|\nabla\vec{L}\cdot \vec{\nu}|\ dx^2=0\quad.
\ee
Let $\vec{w}_\ep:=\phi_\ep\ \vec{\nu}$. We have, using (\ref{II.43})
\be
\label{II.44}
\begin{array}{l}
\ds\p_{\vec{w}_\ep}\lf(\Pi^0_{lk}\circ h\rg)(\vec{\Phi})=-i\,\int_{D^2_\ep}\vec{\nu}\cdot\nabla\vec{L}\,\nabla\phi_\ep\  (A^l)^t\,S_1\, X^k_0\ dx^2\\[5mm]
\ds\quad\quad\quad\quad\quad\quad\quad\quad\quad-i\,\int_{D^2}\vec{\nu}\cdot\nabla^\perp\vec{L}\,\nabla\phi_\ep\   (A^l)^t\,S_2\, X^k_0\ dx^2=O(\ep)\quad,
\end{array}
\ee
and similarly we have
\be
\label{II.45}
\p_{\vec{w}_\ep}\lf(\Pi^0_{lk}\circ h\rg)(\vec{\Phi})=O(\ep)\quad.
\ee
Observe also that $v_\ep=\nabla\vec{w_\ep}\cdot(\nabla\vec{L}+i\nabla^\perp\vec{L})=4\ \p_{\ov{z}}{\phi}_\ep\ \vec{\nu}\cdot\p_z\vec{L}$ converges
strongly to 0 in $L^2$ norm due to (\ref{II.43}) again. Hence we have in particular
\be
\label{II.46}
\lim_{\ep\rightarrow 0} \int\|\nabla c^k_{\vec{w}_\ep}\|^2+|v_\ep\Delta^{-1}\nabla^2 v_\ep|=0\quad.
\ee
We have also
\be
\label{II.47}
\lf|\int_{D^2} \lf[\ov{\nabla}\vec{L}\cdot\nabla\vec{w}_\ep\ \nabla\vec{L}\cdot\nabla\vec{w}_\ep\rg]\ dx^2\rg|\le\, C\ \ep^{-2}\int_{D^2_\ep} |\nabla\vec{L}\cdot\nu|^2\ dx^2=o(1)\quad,
\ee
and similarly
\be
\label{II.48}
\lf|\int_{D^2} \lf[\ov{\nabla}\vec{L}\cdot\nabla\vec{w}_\ep\ \nabla^\perp\vec{L}\cdot\nabla\vec{w}_\ep\rg]\ dx^2\rg|\le\, C\ \ep^{-2}\int_{D^2_\ep} |\nabla\vec{L}\cdot\nu|^2\ dx^2=o(1)\quad,
\ee
Hence, combining (\ref{II.46}), (\ref{II.47}) and (\ref{II.48}) we obtain that
\be
\label{II.49}
\begin{array}{l}
\ds\p^2_{\vec{w}^2_\ep}\lf(\Pi^0_{lk}\circ h\rg)(\vec{\Phi})=
\ds-2^{-1}i\int_{D^2}  |\nabla \vec{L}|^2\ \nabla\phi_\ep\ov{\nabla}\phi_\ep (A^l)^t\, S_1\, X^k_0\ dx^2\\[5mm]
\ds-i\int_{D^2}|\nabla\vec{L}|^2\,\p_{x_1}\phi_\ep\cdot\p_{x_2}\phi_\ep\ (A^l)^t\, S_2\, X^k_0\ dx^2+o(1)
\end{array}
\ee
Since from Wente theorem the conformal factor $|\nabla\vec{L}|^2=2\, e^{-2\la}$ is continuous we have
\be
\label{II.50}
\begin{array}{l}
\ds\p^2_{\vec{w}^2_\ep}\lf(\Pi^0_{lk}\circ h\rg)(\vec{\Phi})=
\ds- i\ e^{-2\la(x^0)}(A^l)^t\, S_1\, X^k_0(x^0) \int_{{\R}^2}  |\p_{x_1}\phi|^2-|\p_{x_2}\phi|^2\ dx^2\\[5mm]
\ds\quad\quad\quad\quad\quad\quad\quad-2\ i\ e^{-2\la(x^0)}(A^l)^t\, S_2\, X^k_0(x^0) \int_{{\R}^2}  \p_{x_1}\phi\,\p_{x_2}\phi\ dx^2+o(1)
\end{array}
\ee
Let $\nu:=2\ e^{2\la(x^0)}\lf(\int_{{\R}^2} \p_{x_i}\phi\p_{x_j}\phi \ dx^2 \rg)dx^i\otimes dx^j$, we have
\be
\label{II.51}
\p^2_{\vec{w}^2_\ep}\lf(\Pi^0_{lk}\circ h\rg)(\vec{\Phi})=-i\, e^{-4\la(x^0)}\ (A^l)^t(x^0)\,\nu^0\, X^k_0(x^0)+o(1)
\ee
Following the computations in previous subsections this gives
\be
\label{II.52}
\lim_{\ep\rightarrow 0}\p^2_{\vec{w}^2_\ep}\lf(\Pi^0_{lk}\circ h\rg)(\vec{\Phi})=-4\,i\, \Im\lf(\om^l\otimes\om^k(x^0),\int_{{\C}}\p\phi\otimes\p\phi\rg)_{wp}\quad.
\ee
And similarly we have
\be
\label{II.53}
\lim_{\ep\rightarrow 0}\p^2_{\vec{w}^2_\ep}\lf(\Pi^1_{lk}\circ h\rg)(\vec{\Phi})=-4\,i \,\Im\lf(i\,\om^l\otimes\om^k(x^0)+\sum_{j=1}^gc_{kj}\,\om^l\otimes\om^j(x^0),\int_{{\C}}\p\phi\otimes\p\phi\rg)_{wp}\quad.
\ee
This closes the proof of proposition~\ref{pr-second}.\hfill$\Box$
\subsubsection{The infinitesimal subjectivity of the 2-jets of linear combinations of $\Pi^0_{kl}\circ h$ and $\Pi^1_{pq}\circ h$ in the corresponding span of  $\p\Pi^0_{kl}$ and $\p\Pi^1_{pq}$ at an isothermic immersion.}
For any vector space $V$ we denote by $S(V\otimes V)$ the subspace of vectors in the tensor product $V\otimes V$ which are symmetric i.e. generated
by $v_1\otimes v_2+v_2\otimes v_1$ for any $(v_1,v_2)\in V^2$.
The goal of this subsection is to prove the following result
\begin{Prop}
\label{pr-surj}
Let $\vec{\Phi}$ be a weak branched immersion of ${\mathcal F}_\Sigma$. Let $(t^j_{kl})$ and $(s^j_{pq})$ for $(k,l)\in \{1\cdots g\}^2$ for $(p,q)\in \{1\cdots g\}^2$ and $j=1\cdots n$ be two families of real numbers such that
the family of $n$ one forms on the space of sections of the bundle $S(T^\ast\Sigma\otimes T^\ast\Sigma)$ given by
\[
L^j(\nu):=\sum_{k,l=1}^gt^j_{kl}\, \p_\nu\Pi^0_{kl}+\sum_{p,q=1}^gs^j_{pq}\,\p_\nu\Pi^1_{pq}\in i\,{\R}
\]
is a free family at $\vec{\Phi}^\ast g_{{\R}^m}$. Assume that the family of corresponding $n$ one forms on $W^{2,2}\cap W^{1,\infty}(\Sigma,{\R}^m)$ given by
\[
{\mathcal L}^j(\vec{w}):=\sum_{k,l=1}^gt^j_{kl}\, \p_{\vec{w}}(\Pi^0_{kl}\circ h)+\sum_{p,q=1}^gs^j_{pq}\,\p_{\vec{w}}(\Pi^1_{pq}\circ h)\in i\,{\R}
\]
at $\vec{\Phi}$ is not free anymore. Then, the rank of the span of $({\mathcal L}^j)_{j=1\cdots g}$ is $n-1$. Moreover introducing the following quadratic forms
on $W^{2,2}\cap W^{1,\infty}(\Sigma,{\R}^m)$ 
\[
{\mathcal Q}^j(\vec{w}):=2^{-1}\sum_{k,l=1}^gt^j_{kl}\, \p^2_{\vec{w}^2}(\Pi^0_{kl}\circ h)+2^{-1}\sum_{p,q=1}^gs^j_{pq}\,\p^2_{\vec{w}^2}(\Pi^1_{pq}\circ h)\in i\,{\R}
\]
there exists a neighborhood $U$ of $0$ in ${\R}^n$ such that, for any $u=(u_1\cdots u_n)\in U$ there exists $\vec{w}\in W^{2,2}\cap W^{1,\infty}(\Sigma,{\R}^m)$ 
\be
\label{II.54}
\forall\, j=1\cdots n\quad\quad i\, u_j={\mathcal L}^j(\vec{w})+{\mathcal Q}^j(\vec{w})\quad.
\ee
\hfill$\Box$
\end{Prop}
\noindent{\bf Proof of proposition~\ref{pr-surj}.}
The first statement is a generalization of proposition~\ref{pr-free} to general families of linear combination of frequencies. Each drop in the rank
of the span of $(L^j)_{j=1\cdots n}$ and the rank of the span of $({\mathcal L}^j)_{j=1\cdots n}$ corresponds to the existence 
of an independent holomorphic quadratic differential $q$ satisfying
\be
\label{II.55}
(q,\vec{\mathfrak h}^0)_{wp}\equiv 0
\ee
and hence, $\vec{\Phi}$ has to be isothermic and, because of proposition~\ref{pr-II.2}, there can be at most one such independent $q$ satisfying (\ref{II.56})
and by consequence almost one drop.
Thus we have
\[
\mbox{rank}({\mathcal L}^j)_{j=1\cdots n}\ge n-1\quad.
\]
If rank $({\mathcal L}^j)_{j=1\cdots n}=n$ this implies that $({\mathcal L}^j)_{j=1\cdots n}$ is a submersion into $i\,{\R}^n$ and ({\ref{II.54}) is clear.

\medskip

We shall assume from now on that $\mbox{rank}({\mathcal L}^j)_{j=1\cdots n}= n-1$ and, without loss of generalities, we can assume that
\be
\label{II.56}
{\mathcal L}^1=0\quad\quad\mbox{ and }\quad\quad\mbox{rank}({\mathcal L}^j)_{j=2\cdots n}=n-1\quad.
\ee
Hence in order to prove (\ref{II.54}) it suffices to find two directions $\vec{w}_\pm$ such that
\be
\label{II.57}
\pm\ i^{-1}\,{\mathcal Q}^1(\vec{w}_\pm)>0\quad\quad\mbox{ and }\quad\quad\forall j=2\cdots n\quad {\mathcal L}^j(\vec{w}_\pm)+{\mathcal Q}^j(\vec{w}_\pm)=0\quad.
\ee
Because of formulas (\ref{II.be26}) and (\ref{II.be27}), since $L^1$ is assumed to be non zero we have that the following holomorphic
quadratic form is non zero
\be
\label{II.58}
q^1:=\sum_{kl}^gt^1_{k,l}\,\om^k\otimes\om^l+\sum_{p,q=1}^gs^1_{pq}\,i\, \om^p\otimes\om^q+\sum_{p,q,r=1}^g s^1_{pq}\,c_{pr}\, \om^p\otimes\om^r\ne 0
\ee
Using proposition~\ref{pr-second} we know that for any function $\phi$ on ${\C}$ there exists a sequence $\vec{w}_\ep$ such that
for any $x^0\in\Sigma$
\be
\label{II.59}
\lim_{\ep\rightarrow 0}{\mathcal L}(\vec{w}_\ep)=0\quad\quad\mbox{and }\quad\quad\lim_{\ep\rightarrow 0}{\mathcal Q}^1(\vec{w}_\ep)=-\,2\,i\ \lf(q^1(x^0),\int_{{\C}}\p\vec{\Phi}\otimes\p\vec{\Phi}\rg)_{wp}\quad.
\ee
Since $q^1$ is non zero $q^1(x^0)$ is non-zero almost everywhere and since the expression  above is invariant under rotations we can assume that 
we are chosing local conformal coordinates such that
$q^1(x^0)=q^1_1(x^0)\ dz^2$ where $q^1_1(x^0)\in{\R}\setminus\{0\}$. Hence we have
\be
\label{II.60}
\lim_{\ep\rightarrow 0}{\mathcal Q}^1(\vec{w}_\ep)=-\,2\,i \,e^{-4\la(x^0)}\ q_1^1(x^0)\ \int_{{\C}}|\p_{x_1}\phi|^2-|\p_{x_2}\phi|^2\ dx^2
\ee
Let $a$ be a compactly supported smooth non zero function on ${\R}$. Choose first $\phi(x_1,x_2):=a(2\,x_1)\,a(x_2)$ we get
\[
\int_{{\C}}|\p_{x_1}\phi|^2-|\p_{x_2}\phi|^2\ dx^2=\frac{3}{2}\,\int_{\R}\dot{a}^2(t)\, dt\ \int_{\R}{a}^2(t)\, dt>0\quad,
\]
whereas by choosing instead $\phi(x_1,x_2):=a(x_1)\,a(2\,x_2)$ we obtain
\[
\int_{{\C}}|\p_{x_1}\phi|^2-|\p_{x_2}\phi|^2\ dx^2=-\frac{3}{2}\,\int_{\R}\dot{a}^2(t)\, dt\ \int_{\R}{a}^2(t)\, dt<0\quad.
\]
Hence, since $a$ is arbitrary in the class of  compactly supported smooth non zero function on ${\R}$, we can find two sequences $\vec{w}_\ep^\pm$ such that
\be
\label{II.61}
\lim_{\ep\rightarrow 0}{\mathcal L}(\vec{w}^\pm_\ep)=0\quad\quad\mbox{and }\quad\quad\lim_{\ep\rightarrow 0}{\mathcal Q}^1(\vec{w}^\pm_\ep)=\pm\, i\quad.
\ee
Let $(\vec{v}_2\cdots\vec{v}_n)$ be $n$ independent elements of $W^{1,\infty}\cap W^{2,2}(\Sigma,{\R}^m)$ such that ${\mathcal L}(\vec{v}_i)$ forms 
a basis of vectors of ${\mathcal L}(W^{1,\infty}\cap W^{2,2}(\Sigma,{\R}^m))=\sum_{j=1}^n{\mathcal L}^j(W^{1,\infty}\cap W^{2,2}(\Sigma,{\R}^m))\ \ep_j=i\,{\R}^{n-1}\subset\,i\,{\R}^n$ where $\ep_j$ denotes the canonical basis of ${\R}^{n}$ .

 For each $j=1\cdots n$ we denote by ${\mathcal B}^j$ the symmetric bilinear form associated to ${\mathcal Q}^j$. Clearly
 \be
 \label{II.62}
 \limsup_{\ep\rightarrow 0}|{\mathcal B}^j(\vec{w}_\ep^\pm,\vec{v})|<+\infty
 \ee
Denote respectively $\hat{{\mathcal L}}:=\sum_{j=2}^n{\mathcal L}^j\ \ep_j$, $\hat{{\mathcal Q}}:=\sum_{j=2}^n{\mathcal Q}^j\ \ep_j$ and  $\hat{{\mathcal Q}}:=\sum_{j=2}^n{\mathcal Q}^j\ \ep_j$. For any $\delta>0$ there exists $\ep_\delta>0$ such that 
\[
\forall \ep<\ep_\delta\quad\quad |{\mathcal L}(\vec{w}_\ep^\pm)|\le\delta\quad.
\]
Consider now the mapping from ${\R}^{n-1}$ into $i\,{\R}^{n-1}$ given by
\[
\begin{array}{l}
\ds\hat{\mathcal N}^\delta_\ep\ :\ (\sigma_2,\cdots,\sigma_n)\longrightarrow \sum_{s=2}^n\sigma_s\,\hat{\mathcal L}(\vec{v}_s)+\delta\sum_{s=2}^n\sigma_s^2\hat{\mathcal Q}( \vec{v}_s)
+\delta\,{\mathcal Q}(\vec{w}_\ep^\pm) +\hat{\mathcal L}(\vec{w}_\ep^\pm)+2\ \delta\ \sum_{s=2}^n\sigma_s\ \hat{\mathcal B}(\vec{w}_\ep^\pm,\vec{v}_s)
\end{array}
\]
Since the linear map 
\[
 (\sigma_2,\cdots,\sigma_n)\longrightarrow \sum_{s=2}^n\sigma_s\,\hat{\mathcal L}(\vec{v}_s)
\]
is an isomorphism, for $\delta$ small enough, $\delta<\delta_0$ and for any $\ep<\ep_\delta$ the local inversion theorem
gives the existence of $(\ov{\sigma}_2,\cdots\ov{\sigma}_n)$ such that
\be
\label{II.63}
\hat{\mathcal N}^\delta_\ep(\ov{\sigma}_2,\cdots\ov{\sigma}_n)=0\quad\quad\mbox{and }\quad\quad(\ov{\sigma}_2,\cdots\ov{\sigma}_n)=O(\delta)
\ee
Hence we have in one hand
\be
\label{II.64}
\forall\ j=2\cdots n\quad\quad{\mathcal L}^j+{\mathcal Q}^j\lf(\delta\sum_{s=2}^n\ov{\sigma}_s\ \vec{v}_s+\delta\,\vec{w}_\ep^\pm\rg)=0
\ee
and in the other hand
\be
\label{II.65}
{\mathcal Q}^1\lf(\delta\sum_{s=2}^n\ov{\sigma}_s\ \vec{v}_s+\delta\,\vec{w}_\ep^\pm\rg)=\pm\ i\ \delta^2\lf(1+O(\delta)\rg)
\ee
The two assertions (\ref{II.64}) and (\ref{II.65}) imply (\ref{II.57}) for $\delta$ small enough, which concludes the proof of proposition~\ref{pr-surj}.\hfill $\Box$

\medskip

\subsection{Writing arbitrary tangent directions in $W^{1,\infty}\cap W^{2,2}(\Sigma,{\R}^m)$  to the pre-image in ${\mathcal E}_\Sigma$ of a sub-manifold of the Teichm\"uller space as a combination of derivatives
of paths of weak immersions within the sub-manifold.}

Combining now proposition~\ref{pr-surj} together with lemma~\ref{lma-2} we are going to obtain theorem~\ref{th-iso-surj} which is one of the main achievement
of the present work.

\begin{Th}
\label{th-iso-surj}
Let ${\mathcal N}$ be a non degenerate smooth ${\R}^n$ valued function on the Teichm\"uller space of $\Sigma$, closed 2-manifold
- i.e. ${\mathcal N}^{-1}\{0\}$ is a submanifold of ${\mathcal T}_\Sigma$. Let $\tau_0$ be a point in
${\mathcal N}^{-1}\{0\}$ such that there exist $n$ independent linear combinations of periods whose derivative are generating $T_\tau({\mathcal N}^{-1}\{0\})$.
Let $\vec{\Phi}$ be a weak immersion of $\Sigma$ into ${\R}^m$ such that $\tau(\vec{\Phi}^\ast g_{{\R}^m})=\tau_0$ and let $\vec{w}$
be an element of $W^{1,\infty}\cap W^{2,2}(\Sigma,{\R}^m)$ tangent to $\tau^{-1}\lf({\mathcal N}^{-1}\{0\}\rg)$ :
\[
\frac{d}{dt}{\mathcal N}(\tau((\vec{\Phi}+t\vec{w})^\ast g_{{\R}^m}))(0)=0\quad.
\]
then there exist two paths 
in $W^{1,\infty}\cap W^{2,2}(\Sigma,{\R}^m)$, $\vec{\Phi}_1(t)$
and $\vec{\Phi}_2(t)$ for $t\in (-1,+1)$, continuous and differentiable at $t=0$ such that
\be
\label{II.65a}
\begin{array}{l}
\ds\vec{\Phi}_1(0)=\vec{\Phi}_2(0)=\vec{\Phi}\quad,\quad\vec{w}=\frac{d}{dt}(\vec{\Phi}_1+\vec{\Phi}_2)(0)\\[5mm]
\ds\quad\mbox{ and }\quad\forall\,t\in(-1,+1)\quad
{\mathcal N}(\tau(\vec{\Phi}_1(t)^\ast g_{{\R}^m}))={\mathcal N}(\tau(\vec{\Phi}_2(t)^\ast g_{{\R}^m}))\equiv 0\quad.
\end{array}
\ee
\hfill $\Box$
\end{Th}
\noindent{\bf Proof of theorem~\ref{th-iso-surj}.}

Let $(t^j_{kl})$ and $(s^j_{pq})$ for $(k,l)\in \{1\cdots g\}^2$ for $(p,q)\in \{1\cdots g\}^2$ and $j=1\cdots n$ be two families of real numbers such that
the family of $n$ one forms on the space of sections of the bundle $S(T^\ast\Sigma\otimes T^\ast\Sigma)$ given by
\[
L^j(\nu):=\sum_{k,l=1}^gt^j_{kl}\, \p_\nu\Pi^0_{kl}+\sum_{p,q=1}^gs^j_{pq}\,\p_\nu\Pi^1_{pq}\in i\,{\R}
\]
is a free family generating the tangent space to ${\mathcal N}^{-1}\{0\}$ : 
\[
\bigcap_{j=1}^nKer (L^j)=T_{h^0}({(\mathcal N\circ\tau)}^{-1}\{0\})\quad.
\]
 Let
\[
{\mathcal L}^j(\vec{w}):=\sum_{k,l=1}^gt^j_{kl}\, \p_{\vec{w}}(\Pi^0_{kl}\circ h)+\sum_{p,q=1}^gs^j_{pq}\,\p_{\vec{w}}(\Pi^1_{pq}\circ h)
\]
be the corresponding family at the level of the immersion $\vec{\Phi}$. If the rank of $({\mathcal L}^j)_{j=1\cdots n}$ is also $n$ then we easily construct
by local inversion theorem a smooth path $\vec{\Phi}(t)$ such that
\[
\vec{\Phi}(0)=\vec{\Phi}\quad\vec{w}=\frac{d}{dt}\vec{\Phi}(0)\quad\mbox{and}\quad{\mathcal N}(\tau(\vec{\Phi}(t)^\ast g_{{\R}^m}))\equiv 0
\]
and the conclusion of the theorem are obtained by taking $\vec{\Phi}_1(t):=\vec{\Phi}(t)$ and $\vec{\Phi}_2(t):\equiv \vec{\Phi}$.

\medskip

Assume now that rank$({\mathcal L}^j)<n$ hence we know from proposition    that  rank$({\mathcal L}^j)=n-1$ and $\vec{\Phi}$ has to be a weak isothermic immersion of ${\mathcal E}_\Sigma$. The given direction $\vec{w}\in W^{2,2}\cap W^{1,\infty}(\Sigma,{\R}^m)$ satisfies
\[
\forall\, j=1\cdots n\quad\quad{\mathcal L}^j(\vec{w}):=\sum_{k,l=1}^gt^j_{kl}\, \p_{\vec{w}}(\Pi^0_{kl}\circ h)+\sum_{p,q=1}^gs^j_{pq}\,\p_{\vec{w}}(\Pi^1_{pq}\circ h)=0
\]
We are then in a position to apply lemma~\ref{lma-2} to the smooth map given by
\[
\vec{\Psi}\in {\mathcal E}_\Sigma\longrightarrow {\mathcal N}\circ\tau\circ h(\vec{\Psi})={\mathcal N}(\tau(\vec{\Psi}^\ast g_{{\R}^m}))\in{\R}^n\quad.
\]
Indeed the Banach manifold ${\mathcal E}_\Sigma$ is in fact an open set of $W^{2,2}\cap W^{1,\infty}(\Sigma,{\R}^m)$.
From the lemma~\ref{lma-2} there exist two elements $\vec{w}_1,\vec{w}_2$ in $W^{2,2}\cap W^{1,\infty}(\Sigma,{\R}^m)$ and two paths $\vec{\Phi}_1(t)$, $\vec{\Phi}_2(t)$
 in $W^{2,2}\cap W^{1,\infty}(\Sigma,{\R}^m)$ which are continuous and differentiable at $0$ such that for $i=1,2$ one has
\be
\label{II.66}
\vec{\Phi}_i(0)=\vec{\Phi}\quad\quad\frac{d\vec{\Phi}_i}{dt}(0)=\vec{w}_i\quad\quad\vec{w}=\vec{w}_1+\vec{w}_2
\ee
and
\be
\label{II.67}
\forall\,t\in(-1,+1)\quad
{\mathcal N}(\tau(\vec{\Phi}_1(t)^\ast g_{{\R}^m}))={\mathcal N}(\tau(\vec{\Phi}_1(t)^\ast g_{{\R}^m}))\equiv 0\quad.
\ee
This concludes the proof of  theorem~\ref{th-iso-surj}.\hfill $\Box$ 
\section{Variations of various energies, Area, Willmore, framed Willmore energies... for weak immersions evolving within
a sub-manifold of the Teichm\"uller Space.}
\subsection{Proof of theorem~\ref{th-0-1}.}
We present the theorem in the case of the Willmore energy only since the proof is identical for any other smooth lagrangian.
It is proved in \cite{Ri3}, \cite{MR2} that $W$ is Frechet differentiable at any point in the Banach manifold ${\mathcal E}_\Sigma$ and that for any $C^1$ path $\vec{\Phi}_1$ 
in ${\mathcal E}_\Sigma$
\be
\label{II.68}
dW(\vec{\Phi})\cdot\vec{w}=\int_{\Sigma}d^{\ast_g}\lf[ d\vec{H}-3\, \pi_{\vec{n}}(d\vec{H})+\star_{{\R}^m}(\ast_{g}d\vec{n}_{\vec{\Phi}}\wedge\vec{H})\rg]\cdot \vec{w}\ dvol_g
\ee
where $d\vec{\Phi}/dt(0)=\vec{w}$. Let $(L^j)_{j=1\cdots n}$ be a free family of linear combination of periods 
\[
L^j(\nu):=\sum_{k,l=1}^gt^j_{kl}\, \p_\nu\Pi^0_{kl}+\sum_{p,q=1}^gs^j_{pq}\,\p_\nu\Pi^1_{pq}\in i\,{\R}
\]
generating the tangent space to $N$ at $\tau(\vec{\Phi})$ : i.e. $T_\tau N=\cap_{j=1}^nKer L^j$. Assume first that the corresponding family of $n$ one-forms on $W^{2,2}\cap W^{1,\infty}(\Sigma)$. Let
\[
{\mathcal L}^j(\vec{w}):=\sum_{k,l=1}^gt^j_{kl}\, \p_{\vec{w}}(\Pi^0_{kl}\circ h)+\sum_{p,q=1}^gs^j_{pq}\,\p_{\vec{w}}(\Pi^1_{pq}\circ h)
\]
be the corresponding family at the level of the immersion $\vec{\Phi}$. If the rank of $({\mathcal L}^j)_{j=1\cdots n}$ is also $n$ then
the constraint is not degenerate and for any $\vec{w}$ satisfying
\be
\label{II.69}
{\mathcal L}^j(\vec{w})=0\quad\forall j=1\cdots n
\ee
we easily construct by local inversion theorem a smooth path $\vec{\Phi}(t)$ such that
$$\vec{\Phi}(0)=\vec{\Phi}\quad\vec{w}=\frac{d}{dt}\vec{\Phi}(0)\quad\mbox{and}\quad \tau(\vec{\Phi}(t))\in N\quad\forall t\quad.$$ Hence, combining (\ref{II.e25}), (\ref{II.f25}) and (\ref{II.68}) we obtain
(\ref{00-2}) and the theorem is proved in the non degenerate case. 

Assuming now that  the rank of $({\mathcal L}^j)_{j=1\cdots n}$ is less than $n$, then we use theorem~\ref{th-iso-surj} and we deduce that for any $\vec{w}$ satisfying (\ref{II.69}) there exist
two paths in $N$, $\vec{\Phi}_1$ and $\vec{\Phi}_2$, differentiable at $t=0$ and satisfying $$\vec{\Phi}_1(0)=\vec{\Phi}_2(0)\quad\mbox{ and }\quad\vec{w}=\frac{d}{dt}(\vec{\Phi}_1+\vec{\Phi}_2)(0)$$. Hence
(\ref{II.68}) holds for such an arbitrary $\vec{w}$ satisfying (\ref{II.69}) and we deduce (\ref{00-2}) in the degenerate case which closes the proof of theorem~\ref{th-0-1}. \hfill $\Box$

\subsection{Proof of theorem~\ref{th-0-2} and theorem~\ref{th-0-4}.}

Once the PDE's (\ref{00-2}) and (\ref{00-3}) are derived, their weak solutions in ${\mathcal F}_\Sigma$ are known to be smooth away from the branched points (see \cite{Ri2}, \cite{Ri3} and \cite{MR3}). Hence theorem~\ref{th-0-2} is proved. Now theorem~\ref{th-0-4} follows from theorem~\ref{th-0-2} and theorem~\ref{th-0-3} (see such an argument in a similar context in \cite{Ri3}
and \cite{MR2} or in \cite{Ri1}). \hfill $\Box$
\appendix
\section{Appendix}
\reset
We shall need the following elementary analysis lemma which is well known but that we recall with a proof for the convenience of the reader.
\begin{Lma}
\label{lma-1}
Let ${\mathcal N}=({\mathcal N}^1,\cdots{\mathcal N}^n)$ be a $C^\infty$ map from a Banach space $E$ into ${\R}^n$.
Denote $\hat{\mathcal N}=({\mathcal N}^2,\cdots{\mathcal N}^n)$ and assume 
\begin{itemize}
 \item[i)] ${\mathcal N}(0)=0$
 \item[ii)] The linear map $\p\hat{\mathcal N}(0)$ on $E$ is a submersion into ${\R}^{n-1}$.
 \item[iii)] The linear form $\p{\mathcal N}^1$ is equal to zero.
 \item[iv)] There exist 2 vectors $\vec{w}_\pm\in E$ such that
 \[
 \p^2_{\vec{w}_{\pm}^2}{\mathcal N}^1(0)=\pm 1\quad\mbox{ and }\quad\quad \p_{\vec{w}_{\pm}}\hat{\mathcal N}(0)=0.
 \]
\end{itemize}
Let $\vec{w}\in E$ such that
\[
\p^2_{\vec{w}^2}{\mathcal N}^1(0)=0\quad\quad\mbox{ and }\quad\quad\p_{\vec{w}}\hat{\mathcal N}(0)=0\quad.
\]
Then there exists a path $\vec{w}(t)\in E$ for $t$ in a neighborhood of zero, continuous and differentiable at the origin, such that
\[
\vec{w}(0)=0\quad,\quad\frac{d\vec{w}}{dt}(0)=\vec{w}\quad\quad\mbox{ and }\quad\forall\, t\quad{\mathcal N}(\vec{w}(t))\equiv 0\quad.
\]
\end{Lma}
\noindent{\bf Proof of lemma~\ref{lma-1}.}
Let $\vec{v}_2\cdots \vec{v}_n$ be family of $n-1$ independent vectors of $E$ such that $(\p_{\vec{v}_s}\hat{\mathcal N}(0))_{s=2\cdots n}$
realizes a basis of the image of $E$ by $\p\hat{\mathcal N}(0)$. Denote by $F$ the sub-vector space in $E$ generated by the $\vec{v}_i$.

We shall look for a path $\vec{w}(t)$ of the form 
\[
\vec{w}(t):=t\,\vec{w}+\sum_{s=2}^n\sigma_s(t)\ \vec{v}_s+\sigma^+(t)\ \vec{w}_++\sigma^-(t)\ \vec{w}_-\quad.
\]
with $\sum_s|\sigma_s(t)|+|\sigma_+(t)|+|\sigma_-(t)|=o(t)$. Denotes by $\hat{\mathcal L}^{-1}$ to be the inverse of $\hat{\mathcal L}$ the restriction
to $F=Span\{\vec{v}_2\cdots\vec{v}_n\}$ into ${\R}^{n-1}$. For $t$, $\sigma^+$ and $\sigma^-$ given and small enough we denote by
$\sigma_s(t,\la_+,\la_-)$ the unique solution obtained by fixed point method of
\[
\sum_{s=2}^{n}\sigma_s(t,\sigma_+,\sigma_-)\ \vec{v}_s=-\hat{\mathcal L}^{-1}\lf(\hat{\mathcal N}\lf(t\,\vec{w}+\sum_{s=2}^n\sigma_s(t,\sigma^+,\sigma^-)\ \vec{v}_s+\sigma^+\ \vec{w}_++\sigma^-\ \vec{w}_-\rg)-\sum_{s=2}^{n}\sigma_s(t,\la_+,\la_-)\ \vec{v}_s\rg)
\]
where we are using that $\p_{\vec{w}_\pm}{\mathcal N}(0)=0$ and $\p_{\vec{w}}{\mathcal N}(0)=0$. The fixed point argument - local inversion theorem -
gives a solution which is smooth w.r.t. $(t,\sigma_+,\sigma_-)$ and satisfying
\be
\label{A-1}
\hat{\mathcal N}\lf(t\,\vec{w}+\sum_{s=2}^n\sigma_s(t,\sigma^+,\sigma^-)\ \vec{v}_s+\sigma^+\ \vec{w}_++\sigma^-\ \vec{w}_-\rg)\equiv 0\quad,
\ee
and
\be
\label{A-2}
|\sigma_s(t,\sigma^+,\sigma^-)|=O(t^2+|\sigma^+|^2+|\sigma^-|^2)\quad.
\ee
It remains to adjust $\sigma^+$ and $\sigma^-$ depending on $t$ in such a way that 
\be
\label{A-3}
{\mathcal N}^1\lf(t\,\vec{w}+\sum_{s=2}^n\sigma_s(t,\sigma^+(t),\sigma^-(t))\ \vec{v}_s+\sigma^+(t)\ \vec{w}_++\sigma^-(t)\ \vec{w}_-\rg)\equiv 0
\ee
and such that 
\be
\label{A-4}
|\sigma^+(t)|+|\sigma^-(t)|=o(t)\quad.
\ee
Let $\vec{S}(t,\sigma^+,\sigma^-):=\sum_{s=2}^n\sigma_s(t,\sigma^+,\sigma^-)\ \vec{v}_s$. Because of (\ref{A-2}) we have
$|\vec{S}(t,\sigma^+,\sigma^-)|=O(t^2+|\sigma^+|^2+|\sigma^-|^2)$.  With the help of this notation we have
\be
\label{A-5}
\begin{array}{l}
\ds{\mathcal N}^1\lf(t\,\vec{w}+\sum_{s=2}^n\sigma_s(t,\sigma^+,\sigma^-)\ \vec{v}_s+\sigma^+\ \vec{w}_++\sigma^-\ \vec{w}_-\rg)\\[5mm]
\ds\quad= {\mathcal Q}^1(t\,\vec{w}+\sigma^+\ \vec{w}_++\sigma^-\ \vec{w}_-)+2\ {\mathcal B}^1\lf(t\,\vec{w}+\sigma^+\ \vec{w}_++\sigma^-\ \vec{w}_-,
\vec{S}(t,\sigma^+,\sigma^-)\rg)+{\mathcal R}^1(t,\sigma^+,\sigma^-)
\end{array}
\ee
where ${\mathcal Q}^1$ is the quadratic form given by ${\mathcal Q}^1:=2^{-1}\p^2{\mathcal N}^1(0)$ and ${\mathcal B}^1$ is the associated bilinear 
form. We have moreover 
\[
|{\mathcal R}^1(t,\sigma^+,\sigma^-)|\le C\, [|t|^3+|\sigma^+|^3+|\sigma^-|^3]\quad.
\]
Consider now 
\[
{\mathcal T}(t,\sigma^+,\sigma^-)=2\ {\mathcal B}^1\lf(t\,\vec{w}+\sigma^+\ \vec{w}_++\sigma^-\ \vec{w}_-,
\vec{S}(t,\sigma^+,\sigma^-)\rg)+{\mathcal R}^1(t,\sigma^+,\sigma^-)\quad.
\]
We have
\[
{\mathcal T}(t,\sigma^+,\sigma^-)={\mathcal T}(t,0,0)+\sum_{\pm}{a}_{\pm}(t)\sigma^\pm+\sum_{\pm}b_{\pm,\pm}(t)\sigma^\pm\sigma^\pm+{\mathcal U}(t,\sigma^+,\sigma^-)
\]
where
\be
\label{A-6}
|{\mathcal T}(t,0,0)|+\sum_{\pm}t\,|{a}_{\pm}(t)|+\sum_{\pm}t^2\,|b_{\pm,\pm}(t)|=O(|t|^3)\quad\mbox{ and }\quad
|{\mathcal U}(t,\sigma^+,\sigma^-)|\le C\ [|\sigma^+|^3+|\sigma^-|^3]
\ee
We have then
\be
\label{A-7}
\begin{array}{l}
{\mathcal N}^1\lf(t\,\vec{w}+\sum_{s=2}^n\sigma_s(t,\sigma^+,\sigma^-)\ \vec{v}_s+\sigma^+\ \vec{w}_++\sigma^-\ \vec{w}_-\rg)
={\mathcal Q}^1(\sigma^+\ \vec{w}_++\sigma^-\ \vec{w}_-)\\[5mm]
\ds\quad+2\sum_{\pm}\ \sigma^\pm\,{\mathcal B}^1(\vec{w},\vec{w}_\pm)
+{\mathcal T}(t,0,0)+2\,\sum_{\pm}{a}_{\pm}(t)\sigma^\pm+\sum_{\pm}b_{\pm,\pm}(t)\sigma^\pm\sigma^\pm+{\mathcal U}(t,\sigma^+,\sigma^-)
\end{array}
\ee
In order to chose $\sigma^\pm(t)$ such that (\ref{A-3}) and (\ref{A-4}) hold we are going to consider two alternatives separartly.

\medskip

\noindent{\bf 1st case :} either ${\mathcal B}^1(\vec{w},\vec{w}_+)\ne 0 $ or ${\mathcal B}^1(\vec{w},\vec{w}_-)\ne 0 $.

\medskip

Assume for instance $a_+^0=2\,{\mathcal B}^1(\vec{w},\vec{w}_+)\ne 0$. 
In this case we choose $\sigma^-(t)\equiv 0$. and in order to ensure (\ref{A-3}) we are looking for $\sigma^+(t)$ satisfying
\[
0=(1+b_{++}(t))(\sigma^+)^2+(t\,a_+^0+2\,a_+(t))\ \sigma^+(t)+{\mathcal T}(t,0,0)+{\mathcal U}(t,\sigma^+,0)\quad,
\]
or in other words
\be
\label{A-8}
\begin{array}{l}
\ds\sigma^+=-(a_+^0+2\,t^{-1} a_+(t))^{-1}\,t^{-1}{\mathcal T}(t,0,0)\\[5mm]
\ds\quad\quad-(a_+^0+2\,t^{-1}\,a_+(t))^{-1}\ \lf[t^{-1}(1+b_{++}(t))(\sigma^+)^2+{\mathcal U}(t,\sigma^+,0)\rg]\quad.
\end{array}
\ee
Since
\[
\lf|(a_+^0+2\,t^{-1} a_+(t))^{-1}\,t^{-1}{\mathcal T}(t,0,0)\rg|=O(t^2)
\]
and since
\[
\lf|(a_+^0+2\,a_+(t))^{-1}\ \lf[t^{-1}(1+b_{++}(t))(\sigma^+)^2+t^{-1}{\mathcal U}(t,\sigma^+,0)\rg]\rg|=O(t^{-1}\,(\sigma^+)^2)
\]
There is a unique solution to (\ref{A-8}) satisfying $\sigma^+(t)=O(t^2)$ and by taking this $\sigma^+(t)$ we just constructed and $\sigma^-(t)\equiv 0$, (\ref{A-3}) and (\ref{A-4}) 
is satisfied and the lemma~\ref{lma-1} is proved in this case.

\medskip

\noindent{\bf 2nd case :}  ${\mathcal B}^1(\vec{w},\vec{w}_+)={\mathcal B}^1(\vec{w},\vec{w}_-)= 0 $. 

\medskip

For any $t$ such that ${\mathcal T}(t,0,0)=0$ we take $\sigma^+(t)=\sigma^-(t)= 0$. For $t$ such that ${\mathcal T}(t,0,0)>0$ we take $\sigma^+(t)=0$ and we look for $\sigma^-$ solution to
\be
\label{A-9a}
0=-\lf(1-b_{--}(t)-(\sigma^-)^{-2}{\mathcal U}(t,0,\sigma^-)\rg)\,(\sigma^-)^2+2\, a_-(t)\ \sigma^-+{\mathcal T}(t,0,0)
\ee
Since $(\sigma^-)^{-2}\,{\mathcal U}(t,0,\sigma^-)=O(\sigma^-)$, a straightforward argument gives a unique solution to (\ref{A-9}) as being $1+o(t)$ times 
\[
\sigma^-(t):=(1+o(t))\ \lf[-a_-(t)+\ep(t)\sqrt{a_-^2(t)+(1-b_{--}(t))\, {\mathcal T}(t,0,0)}\rg]\lf(1-b_{--}(t)\rg)^{-1}=O(|t|^{3/2})
\]
where $\ep(t)=sign (a_-(t))$ if $a_-(t)\ne 0$ or $\ep(t)=1$ otherwise, which is one of the 2 solutions of
\be
\label{A-9}
0=-\lf(1-b_{--}(t)\rg)\,(\sigma^-)^2+2\, a_-(t)\ \sigma^-+{\mathcal T}(t,0,0)\quad.
\ee
If ${\mathcal T}(t,0,0)<0$ we exchange the roles of $\sigma^-$ and $\sigma^+$. Observe that the path $\vec{w}(t)$ constructed in this way can have discontinuities at the non zero $t$ such that $a_\pm(t)$ vanish. However it fulfills all the conclusions of lemma~\ref{lma-1} in this second and last case
\hfill $\Box$

\medskip

Form lemma~\ref{lma-1} we deduce the following
\begin{Lma}
\label{lma-2}
Let ${\mathcal N}=({\mathcal N}^1,\cdots{\mathcal N}^n)$ be a $C^\infty$ map from a Banach space $E$ into ${\R}^n$.
Denote $\hat{\mathcal N}=({\mathcal N}^2,\cdots{\mathcal N}^n)$ and assume 
\begin{itemize}
 \item[i)] ${\mathcal N}(0)=0$
 \item[ii)] The linear map $\p\hat{\mathcal N}(0)$ on $E$ is a submersion into ${\R}^{n-1}$.
 \item[iii)] The linear form $\p{\mathcal N}^1$ is equal to zero.
 \item[iv)] There exist 2 vectors $\vec{w}_\pm\in E$ such that
 \[
 \p^2_{\vec{w}_{\pm}^2}{\mathcal N}^1(0)=\pm 1\quad\mbox{ and }\quad\quad \p_{\vec{w}_{\pm}}\hat{\mathcal N}(0)=0.
 \]
\end{itemize}
Let $\vec{w}\in E$ such that
\[
\p_{\vec{w}}\hat{\mathcal N}(0)=0\quad.
\]
Then there exists two vectors $\vec{w}_1$ and $\vec{w}_2$ such that 
\be
\label{A-10}
 \p_{\vec{w}_1}\hat{\mathcal N}(0)=\p_{\vec{w}_2}\hat{\mathcal N}(0)=0\quad,\quad\p^2_{\vec{w}_1}{\mathcal N}^1(0)=\p^2_{\vec{w}_1}{\mathcal N}^1(0)=0\quad
,\quad\mbox{ and }\quad\vec{w}=\vec{w}_1+\vec{w}_2
\ee
furthermore there exists 2 paths $\vec{w}_i(t)\in E$ for $t$ in a neighborhood of zero, continuous and differentiable at the origin, such that for $i=1,2$
\[
\vec{w}_i(0)=0\quad,\quad\frac{d\vec{w}_i}{dt}(0)=\vec{w}_i\quad\quad\mbox{ and }\quad\forall\, t\quad{\mathcal N}(\vec{w}_i(t))\equiv 0\quad.
\]
\end{Lma}
\noindent{\bf Proof of lemma~\ref{lma-2}.}
Assume $\p^2_{\vec{w}^2}{\mathcal N}^1(0)=0$ the we take $\vec{w}_1=\vec{w}$ and we apply lemma~\ref{lma-1} to this vector. $\vec{w}_2(t)$ is then 
chosen to be the trivial path $\vec{w}_2(t)\equiv 0$.

\medskip

Assume now $\p^2_{\vec{w}^2}{\mathcal N}^1(0)\ne 0$ and consider for instance the case $\p^2_{\vec{w}^2}{\mathcal N}^1(0)>0$. Consider then the 
2-dimensional vector space given by $Span\{\vec{w},\vec{w}_-\}\subset Ker(\p\hat{\mathcal N})$ - the two vectors $\vec{w}$ and $\vec{w}_-$ cannot be parallel to each other since 
$\p^2_{\vec{w}^2}{\mathcal N}^1(0)$ and $\p^2_{\vec{w}_-^2}{\mathcal N}^1(0)$ have opposite sign. The quadratic form on this $2-$plane 
given by $\p^2{\mathcal N}^1(0)$ has signature $(+,-)$ and therefore $Span\{\vec{w},\vec{w}_-\}$ is generated by a basis of two isotropic vectors
$\vec{e}_1$ and $\vec{e}_2$ satisfying then $\p^2_{\vec{e}_i^2}{\mathcal N}^1(0)=0$. There exist $w_1, w_2\in{\R}$ such that $\vec{w}=w_1\,\vec{e}_1+w_2\,\vec{e}_2$. Then $\vec{w}_1:=w_1\,\vec{e}_1$ and $\vec{w}_2:=w_2\,\vec{e}_2$ satisfy (\ref{A-10}). We then apply lemma~\ref{lma-1} to each of the two
$\vec{w}_i$ and lemma~\ref{lma-2} is proved in this case. 

The last case $\p^2_{\vec{w}^2}{\mathcal N}^1(0)<0$ is solved identically to the previous one after having replaced $\vec{w}_-$ by $\vec{w}_+$.
This concludes the proof of the lemma~\ref{lma-2}.\hfill $\Box$

\end{document}